\documentclass[11pt]{amsart}

\usepackage{amssymb}
\usepackage{latexsym}
\usepackage{amsmath}
\usepackage{amsthm}
\usepackage{amsfonts}
\usepackage{color}
\usepackage{pdfsync}
\usepackage[usenames,dvipsnames]{pstricks}
\usepackage{epsfig}
\usepackage{pst-grad} 
\usepackage{pst-plot} 

\newcommand{\R}{\mathbb{R}}
\newcommand{\N}{\mathbb{N}}

\newcommand{\I}{\mathcal{I}}

\theoremstyle{plain}
\newtheorem{defi}{Definition}[section]
\newtheorem{prop}[defi]{Proposition}
\newtheorem{teo}[defi]{Theorem}
\newtheorem{cor}[defi]{Corollary}
\newtheorem{lema}[defi]{Lemma}

\newtheorem{remark}[defi]{Remark}

\theoremstyle{definition}

\theoremstyle{remark}

\numberwithin{equation}{section}

\begin{document}

\title[]{Lipschitz Regularity for Censored Subdiffusive Integro-Differential Equations with Superfractional Gradient Terms.}

\author[]{Guy Barles}
\address{
Guy Barles:
Laboratoire de Math\'ematiques et Physique Th\'eorique (UMR CNRS 7350), F\'ed\'eration Denis Poisson (FR CNRS 2964),
Universit\'e Fran\c{c}ois Rabelais Tours, Parc de Grandmont, 37200 Tours, France. {\tt Guy.Barles@lmpt.univ-tours.fr} }

\author[]{Erwin Topp}
\address{
Erwin Topp:
Departamento de Ingenier\'\i a Matem\'atica (UMI 2807 CNRS), Universidad de Chile, Casilla 170, Correo 3, Santiago, Chile.
{\tt etopp@dim.uchile.cl} }

\date{\today}

\begin{abstract} 
In this paper we are interested in integro-differential elliptic and parabolic equations involving nonlocal operators 
with order less than one, and a gradient term whose coercivity growth makes it the leading term in the equation.
We obtain Lipschitz regularity results for the associated stationary Dirichlet problem in the case when the nonlocality 
of the operator is confined to the domain, feature which is known in the literature as censored nonlocality. 
As an application of this result, we obtain strong comparison principles which allow us to prove the well-posedness of both the stationary and evolution problems, and steady/ergodic large time behavior for the associated evolution problem.
\end{abstract}


\maketitle


\section{Introduction.}

In \cite{BKLT}, the authors, in collaboration with O. Ley and S. Koike, investigate regularity properties 
for subsolutions of integro-differential stationary equations with a super-fractional gradient term, providing 
analogous results for nonlocal equations to the ones of Capuzzo-Dolcetta, Porretta and Leoni~\cite{Capuzzo-Dolcetta-Leoni-Porretta} 
for superquadratic degenerate elliptic second-order pdes (see also Barles \cite{Barles1}). 
In the present work, our aim is to obtain analogous results but for integro-differential operators 
of order $\sigma <1$, still with a super-fractional gradient term, but in the (intriguing) case of
\textsl{censored operators} set in bounded domains.

In order to be more specific, we consider the model problem
\begin{equation}\label{eq}
\lambda u(x) + (-\Delta)^{\sigma/2}_c u(x) + b(x)|Du(x)|^m = f(x) \quad \mbox{in} \ \Omega,
\end{equation}
where $\Omega \subset \R^N$ is a bounded domain, $\lambda \geq 0$ and $b, f: \bar{\Omega} \to \R$ are
continuous functions with $b(x)>0$ on $\bar{\Omega}$. For $\sigma \in (0,2)$, 
the integro-differential operator $(-\Delta)^{\sigma/2}_c$ is 
known in the literature as the \textsl{censored fractional Laplacian of order $\sigma$}, 
and is defined through the expression
\begin{equation}\label{censoredLaplacian}
(-\Delta)_c^{\sigma/2}u(x) = C_{N, \sigma} \ \mathrm{P.V.} \int_{x + z \in \Omega} [u(x + z) - u(x)] |z|^{-(N + \sigma)} dz,
\end{equation}
where $\mathrm{P.V.}$ stands for the Cauchy principal value and $C_{N, \sigma} > 0$ is a well-known 
normalizing constant, see~\cite{BBC, GQY2}. The ``censored'' appellative is referred to the fact 
that the integration on the set $\{x + z \in \Omega\}$ makes the jumps outside $\Omega$ being indeed censored. 

In the scope of this paper, besides the mentioned censored nonlocality feature of the problem, the 
main assumptions are $0< \sigma < 1$ (subdiffusive operator of order less than $1$), 
and $m > \sigma$ (superfractional coercivity condition). The study of this case is motivated by two principal reasons: first, 
under such conditions we are able to obtain regularity properties that are (maybe surprisingly) 
more sophisticated than in \cite{BKLT}; indeed, we are not only able to obtain global H\"older continuity but also 
global Lipschitz regularity for bounded subsolutions of (\ref{eq}). Secondly, we obtain  comparison principle and
well-posedness of the Dirichlet problem both for stationary and evolution equations, namely
\begin{equation}\label{ev-eq}
u_t + \lambda u + (-\Delta)^{\sigma/2}_c u + b(x) \ |Du|^m = f(x) \quad \mbox{in} \ \Omega \times (0, \infty).
\end{equation}

Concerning Dirichlet problems for nonlocal equations, we remark that in general 
the \textsl{Dirichlet boundary condition} has to be imposed on the complementary of $\Omega$. 
Typically, if we replace $(-\Delta)_c^{\sigma/2}$ in~\eqref{eq} or \eqref{ev-eq} by the fractional Laplacian $(-\Delta)^{\sigma/2}$ defined as
\begin{equation*}
(-\Delta)^{\sigma/2}u(x) = C_{N, \sigma} \ \mathrm{P.V.} \int_{\R^N} [u(x + z) - u(x)] |z|^{-(N + \sigma)} dz, 
\end{equation*}
then this requires the value of the function in the whole space to be evaluated, and the Dirichlet problem reads
\begin{equation}\label{noncensored}
\left \{ \begin{array}{rll} \lambda u(x) + (-\Delta)^{\sigma/2} u(x) + b(x)|Du|^m & = f(x) \quad & \mbox{in} \ \Omega, \\
u & = \varphi \quad & \mbox{in} \ \Omega^c .\end{array} \right .
\end{equation}

On the contrary, in the censored case, since we use only the values of $u$ in $\bar{\Omega}$, we can complement~\eqref{eq} with a \textsl{classical} boundary condition
\begin{equation}\label{bc}
u = \varphi \quad \mbox{on} \ \partial \Omega,
\end{equation}
for a boundary data $\varphi \in C(\partial \Omega)$. 

A third interesting question that can be handled with the development of regularity and comparison 
properties on the current setting 
is the study of the large time behavior for Cauchy-Dirichlet problem associated to~\eqref{ev-eq}. 
We both study the cases when the censored parabolic problem has a steady state asymptotic behavior, and the case
 of the \textsl{ergodic large time behavior}, situation in which we have to solve a stationary problem with 
state-constraint boundary condition (the \textsl{ergodic problem}).

We want to mention immediately a very important point related to our regularity and well-posedness issues
(and, as a consequence, the large time behavior issues) and which justifies our choice of the 
parameters $0 < \sigma < 1$ and $m > \sigma$ : 
under these conditions,
we are able to solve the censored Dirichlet problem in its full generality, but this is because we use in a 
key way the regularity result for~\eqref{eq}. 
On the contrary, in the case $\sigma \geq 1$ and $m > \sigma$, we are unable to prove that 
the censored Dirichlet problem is well-posed, even with the regularity results of \cite{BKLT}. 

Now we detail the different points discussed above, starting with regularity. In this respect, we
remark the results given in~\cite{Capuzzo-Dolcetta-Leoni-Porretta, Barles1} 
are concerned with superquadratic second-order degenerate elliptic problems like
\begin{equation}\label{superquadratic}
\lambda u - \mathrm{Tr}(A (x)D^2 u) + b(x) |Du|^m = f(x) \quad \mbox{in}  \ \Omega,
\end{equation}
i.e. when $m>2$. In~\cite{Capuzzo-Dolcetta-Leoni-Porretta, Barles1}, the authors prove that 
if $u : \Omega \to \R$ is a bounded \textsl{viscosity subsolution} of (\ref{superquadratic}), 
then $u$ is locally H\"older continuous with exponent $\alpha := (m-2)(m-1)^{-1}$ and the local H\"older 
seminorm depends only on the data
($L^\infty$ bounds on $A,b$ and upper bound on $f$) and $||\lambda u^-||_\infty$, where $u^-=\min(u,0)$. 

In many interesting situations ($\lambda = 0$ or cases where we have a bound on $\lambda u$ which depends only on the data), 
this H\"older seminorm 
does not really depend on any $L^\infty$ bound nor oscillation of $u$ and 
actually provides an estimate on the $L^\infty$ norm of $u$. The H\"older exponent $(m - 2)(m - 1)^{-1}$ just comes from a simple balance of powers in~\eqref{superquadratic} and this H\"older regularity can be extended up to the boundary of the domain if it is regular enough.

This strategy was recently applied in~\cite{BKLT} 
to get regularity results to integro-differential problems for which~\eqref{eq} is a particular case. However, the presence
of the nonlocal term has an effect on the results in two main directions:
first, the global H\"older exponent found by this method does not follow anymore the ``natural'' balance of 
powers in~\eqref{eq}. Related to this, we must take into account 
the fact that (super)solutions to superfractional nonlocal problems (censored or not)
may not satisfy the boundary condition in the classical sense (we come back later to this fact in a more detailed way). 
This creates a difficulty which is more evident in the non-censored setting because in this framework
we must consider the evaluation of the nonlocal term on functions 
developing boundary jump discontinuities, and this gives way to unbounded terms near the boundary. 
A second consequence of the non locality of the operator is that the H\"older regularity results of~\cite{BKLT} 
do not provide, in general, an estimate on the $L^\infty$ norm of $u$ but on the contrary rely on this $L^\infty$ norm. 
This is where censored problems come into play: in these cases, and if $m>1$, the H\"older regularity 
results of~\cite{BKLT} provide an estimate on the $L^\infty$ norm of $u$ (or more precisely on its oscillation).

In this article, we obtain the global H\"older regularity results of~\cite{BKLT} in 
the case $\sigma <1$ and $m>\sigma$ through a simpler proof, and we extend these results
to global Lipschitz continuity for bounded subsolutions to~\eqref{eq}. 
As in~\cite{BKLT}, the censorship of the operator leads us to a control 
of the oscillation for subsolutions to~\eqref{eq} when $m > 1$.
Our regularity results strongly rely on the leading effect of the gradient 
terms and therefore cannot be applied to the ``critical'' case $m = \sigma$. We refer to~\cite{Silvestre, Chang-Lara-Davila} for regularity results in the case $\sigma = m = 1$ : their arguments are based on the (uniform)
ellipticity of the nonlocal operator.  Finally, the case $m = \sigma < 1$ seems to be difficult since the ellipticity of the operator is not strong enough.

The next step is to consider the well-posedness of the (stationary and evolution) nonlocal Dirichlet problem, 
and as we already mentioned above, 
we must deal with the presence of \textsl{loss of the boundary condition}.
We shall mention that this phenomena also arises in local pdes like~\eqref{superquadratic} because of the 
leading effect of the gradient term.
Indeed, the \textsl{classical} Dirichlet problem (with boundary data being really satisfied by the solution) cannot be solved 
in general and one has to use the \textsl{generalized} Dirichlet problem in the sense of viscosity solutions. We refer to the 
Users' guide~\cite{usersguide} and references therein for an introduction of this concept and to \cite{Barles-DaLio, Lasry-Lions} 
for the applications to~\eqref{superquadratic}-\eqref{bc} where it is shown that the generalized Dirichlet problem is well-posed 
in $C(\bar{\Omega})$ or $C(\bar{\Omega}\times [0,T])$, and with examples where the solution is different from $\varphi$ at points 
on $\partial \Omega$.

For nonlocal equations, the possible loss of the boundary conditions was first studied systematically in Barles, 
Chasseigne and Imbert~\cite{Barles-Chasseigne-Imbert} where an explicit example of such loss of boundary condition 
is provided for a non-censored operator. The first study of Dirichlet problems for nonlocal equations with loss 
of boundary conditions was done in~\cite{Topp}, where well-posedness for non-censored nonlocal Hamilton-Jacobi problems 
which are not necessarily coercive in the gradient are obtained, and then in~\cite{BT} in the non-censored but coercive case. 
As it can be seen in~\cite{BT, Topp}, the non-censorship of the operator allows to incorporate the exterior data $\varphi$ 
in~\eqref{noncensored} into the equation. Roughly speaking, this procedure creates a ``new''  term coming 
from the integration outside the domain which can be regarded 
as an extra discount factor and that turns the problem into a strictly proper one, from which the solvability of the Dirichlet 
problem can be obtained even if the $\lambda$ term is negative.

The current censored setting shares some aspects of the pde framework since the problem is really set on $\bar{\Omega}$ and 
no extra information outside $\bar{\Omega}$ comes into play. However, the censorship of the operator must be regarded 
as an $x$-dependence, and it is known that such dependence creates special difficulties 
not only in the study of well-posedness~\cite{Barles-Imbert}, but also in  regularity~\cite{Barles-Chasseigne-Ciomaga-Imbert-lip}. 
It is because of this difficulty that the case $\sigma \geq 1$ is not treated in the full generality neither for Dirichlet nor 
for Neumann problems (see~\cite{bcgj}). Of course, in the Dirichlet case, the difficulty comes from the loss of boundary conditions, otherwise a more standard comparison result is enough to prove the well-posedness : this is for example the case for  $\sigma \geq 1, m\leq \sigma$.

In this article the Lipschitz continuity of subsolutions allows us to turn around these difficulties and to obtain a strong comparison principle for~\eqref{eq}-\eqref{bc} and its evolution counterpart. In this last case, we have to use a (classical) regularization in time to reduce to a situation where we have Lipschitz continuity in $x$. The key point is that this regularity permits to control the estimates of the nonlocal term when it is of order $\sigma < 1$. 

We remark that the global H\"older regularity results of~\cite{BKLT} would not be enough to treat 
the case $\sigma \geq 1$ and to treat the case $\sigma < 1$ would have been more involved.
Moreover, we recall that the generalized notion of solution involves an evaluation of the equation on the
boundary, and such an evaluation has no clear meaning in superdiffusive problems (that is when $\sigma \geq 1$ 
in~\eqref{censoredLaplacian}) mainly because of the asymmetry of the domain of integration. In fact, some special 
requirements on the test functions must be considered in this context, see~\cite{bcgj, GQY2}.

We finish with the application of regularity and comparison results in the study of large time behavior 
of the solutions of the evolution problem.
In the case when $\lambda > 0$, this problem has the expected convergence to the steady-state solution, while if $\lambda = 0$ and $m > 1$, we prove that the asymptotic behavior of this problem resembles the second-order case studied by T. Tchamba in~\cite{Tchamba}, i.e. a suitable ergodic asymptotic behavior
$$ u(x,t) = -ct + u_\infty (x) + o(1)\; ,$$
where  the function $u_\infty$ is a solution to the ergodic problem.

The steady-state asymptotic behavior is closely related to the uniqueness of the associated stationary problem 
and the strict positiveness of $\lambda$ is the nondegeneracy condition for this problem. 
It is interesting to compare this steady-state behavior with the corresponding asymptotic behavior 
in non-censored problems. As we mentioned above, the non-censorship of the operator implies the existence of an extra 
discount factor leading to the mentioned nondegeneracy condition even if $\lambda$ is negative, 
provided the nonlocal operator recovers sufficient information from $\Omega^c$.

A qualitatively different asymptotic behavior is observed in the case $\lambda = 0$ and $m > 1$, case in which 
we obtain robust regularity results leading to the solvability of the ergodic problem similarly to the periodic setting, 
where no boundary is considered, see~\cite{Barles-Chasseigne-Ciomaga-Imbert, BKLT}. However, here we face an ergodic problem 
which is a \textsl{state-constraint} problem, as in ~\cite{Tchamba}.
The study of the ergodic problem is closely connected to the study of the (stationary) Dirichlet problem, while, 
again as in ~\cite{Tchamba}, the convergence of $u(x,t)+ct$ to $u_\infty$ relies on a Strong Maximum Principle type argument. 
The Strong Maximum Principle we use is inspired by Coville~\cite{coville} (see also Ciomaga~\cite{ciomaga}), and it is 
more related to a topological property of the the support of the measure defining the nonlocal operator than with its ellipticity.

The article is organized along the same lines as this introduction : in Section~\ref{defsolsection}
we present the notion of solution we use in the paper. In Section~\ref{regsection}, we provide 
our regularity results, proving first local Lipschitz continuity, then global H\"older regularity results and 
finally a global Lipschitz continuity result. Section~\ref{cpw} is devoted to state and prove the main Strong 
Comparison Result and to deduce the well-posedness of the initial boundary value problem. In Section~\ref{sec:ep}, we study the ergodic problem, obtaining the existence and uniqueness (up to a constant) of the solutions invoking the Strong Maximum Principle. Finally we describe the different large time behavior of the solutions of the initial boundary value problem in Section~\ref{sec:ltb}.

\medskip

\noindent
{\bf \textit{Basic Notation.}} 
Throughout this paper we denote $Q = \Omega \times (0,+\infty)$ and $\partial Q = \partial \Omega \times (0,+\infty)$. 

For $x \in \R^N$ and $r > 0$, we denote $B_r(x)$ the open ball centered at $x$ with 
radius $r$, $B_r$ if $x = 0$, and $B$ if additionally $r = 1$.

For $x \in \mathcal{O}$, we denote $d_\mathcal{O}(x) = dist(x, \partial \mathcal{O})$. In the case $\mathcal{O} = \Omega$
we simply write $d(x)$. For $\delta > 0$, we denote
$\mathcal{O}_\delta = \{ x \in \mathcal{O} : d_{\mathcal{O}}(x) < \delta \}$. 
In the sequel, we always assume that the boundary of the domain is $C^1$ which means we can find $\delta_0 > 0$ such that the distance function is of class $C^1$ in $\Omega_{\delta_0}$, see~\cite{G-T}.

Finally, we write $u.s.c.$ to refer upper semicontinuous functions and in the same way for $l.s.c.$.


\section{Notion of Solution.}
\label{defsolsection}

For the sake of completeness of the paper, we briefly present the notion of viscosity solution for the ``generalized Dirichlet problem'' in the case of equations like \eqref{ev-eq}. Because of the ``censoring'' in the nonlocal term, this notion is a straightforward extension of its second-order counterpart (See \cite{usersguide}), contrarily 
to the case of general nonlocal equations where the Dirichlet boundary conditions has to be imposed on the complementary of the domain and therefore the definition has to be modified accordingly (See for example \cite{Barles-Chasseigne-Imbert}).

For each $x \in \bar{\Omega}$, we consider $\nu_x$ a nonnegative regular measure of $\R^N$.
The nonlocal operator has the general form
\begin{equation}\label{operator}
\I(u, x) = \int_{\R^N} [u(x + z) - u(x)] \nu_x(dz),
\end{equation}
for a function $u: \R^N \to \R$ for which the above expression has a sense.

We start with the precise basic assumptions on the nonlocal term. The first one is concerned with the ``censoring'' property.

\medskip
\noindent
{\bf (M0-1)} \textsl{If
$
\Omega_\nu =  \bigcup_{x \in \bar{\Omega}} \{ x + \mathrm{supp} \{ \nu_x \} \},
$
}
\textsl{the family $\{ \nu_x \}_x$ satisfies
\begin{equation*}
\Omega_\nu \subseteq \bar{\Omega}.
\end{equation*}
}
Assumption (M0-1) means that, for each $x \in \bar{\Omega}$, $x+z \in \bar{\Omega}$ if $z\in \mathrm{supp} \{ \nu_x \}$
and then the domain of integration in~\eqref{operator} can be replaced just by $\bar{\Omega} - x$. For this reason 
we say the nonlocal operator $\I$ has a \textsl{censored nature} since the jumps outside  $\bar{\Omega}$ are actually censored.

The second key assumption describes the set of functions for which the nonlocal term makes sense and its continuity properties in $x$.

\medskip
\noindent
{\bf (M0-2)} \textsl{If $u \in L^\infty(\bar \Omega) \cap C^1(B_r(x))$ for some $r > 0$, $\I(u, x)$ is well-defined. Moreover, for such $u$ and for any $D \subset \R^N$
\begin{equation*}
\begin{split}
\I[D](u, y) & = \int_{D} [u(y + z) - u(y)] \nu_x(dz), \\
\end{split}
\end{equation*}
is continuous in $B_r(x)$.
}

(M0-2) contains two important informations : on one hand, $I(u,x)$ is an operator which is defined, in particular, on $C^1(\bar{\Omega})$ and therefore it can be seen as an operator of order less than $1$. And if $u$ is $C^1$ then $I(u,x)$ is continuous in $x$, a key property to define viscosity solutions.

In the sequel, we always assume these two properties to be satisfied and we denote by {\bf (M0)} the two assumptions (M0-1) and (M0-2).

Next we consider a parabolic problem with general form
\begin{equation}\label{eqdefi}
\left \{ \begin{array}{rll} u_t - \I(u(\cdot, t), x) + H(x, u, Du) &= 0, \ &\mbox{in} \ Q_T \\
u &= \varphi, \ &\mbox{in} \ \partial Q_T \\
u &= u_0, \ &\mbox{in} \ \bar{\Omega} \times \{ 0 \}. \end{array} \right .
\end{equation}
for continuous data $H, \varphi$ and $u_0$. For $T > 0$, we denote $Q_T = \Omega \times (0,T)$
and $\partial Q_T = \partial \Omega \times (0, T]$. In order to slightly simplify the definition below and since this assumption plays a key role in the existence and uniqueness of a {\em continuous solution},  we also request the following \textsl{compatibility condition} between $\varphi$ and $u_0$
\begin{eqnarray}\label{compatibility}
u_0 (\cdot) \equiv \varphi(\cdot, 0) \quad \mbox{on} \ \partial \Omega.
\end{eqnarray}

For $(x,t) \in \bar{Q}_T$, 
$\delta > 0$, and functions $u, \phi: \bar{Q}_T \to \R$ we write
\begin{equation*}
\begin{split}
\mathcal{E}_\delta(u, \phi, x, t) = & \ \phi_t(x,t) - \I[B_\delta](\phi(\cdot, t), x) - \I[B_\delta^c](u(\cdot, t), x) \\
& + H(x, u(x,t), D\phi(x,t)),
\end{split}
\end{equation*}
whenever this expression has a sense.

Now we present the notion of solution for our parabolic problem~\eqref{eqdefi}. This definition 
can be readily extended to the infinite horizon setting (that is $T = +\infty$) and for the corresponding stationary problem.
\begin{defi}\label{defsol}
An u.s.c. function $u: \bar Q_T \to \R$ is a viscosity subsolution to~\eqref{eqdefi} if for each $(x_0, t_0) \in \bar Q_T$
and each $\phi \in C^1(\bar Q_T)$ such that $(x_0, t_0)$ is a maximum point for $u - \phi$ on $\bar Q_T$, 
for all $\delta > 0$ we have
\begin{equation*}
\left \{ \begin{array}{rll} \mathcal{E}_\delta(u, \phi, x_0, t_0) \leq & 0 \quad & \mbox{if} \ (x_0, t_0) \in Q_T, \\
\max \{ \mathcal{E}_\delta(u, \phi, x_0, t_0), u(x_0, t_0) - \varphi(x_0, t_0) \} \leq & 0 \quad & \mbox{if} \ 
(x_0, t_0) \in \partial Q_T, \\
\max \{ \mathcal{E}_\delta(u, \phi, x_0, t_0), u(x_0, t_0) - u_0(x_0) \} \leq & 0 \quad & \mbox{if} \ t_0 = 0.
\end{array} \right .
\end{equation*}

In the analogous way we define supersolutions and solutions to~\eqref{eqdefi}.
\end{defi}

For reasons that will be made clear later, we also require a definition of viscosity solution for the following
\textsl{state constraint problem}
\begin{eqnarray*}
\left \{ \begin{array}{rll} \partial_t u - \I(u, x) + H(x,u,Du) & = 0 \quad & \mbox{in} \ Q_T, \\
\partial_t u - \I(u, x) + H(x,u,Du) & \geq 0 \quad & \mbox{in} \ \partial Q_T, \\
u & = u_0 \quad & \mbox{on} \ \bar{\Omega} \times \{ 0 \}. \end{array} \right .
\end{eqnarray*}

For a solution of this problem we mean a function $u$ which is a viscosity subsolution to the equation in $Q_T$, 
a viscosity supersolution to the equation in $\bar Q_T$, and which satisfies the initial condition in the generalized
sense presented in Definition~\ref{defsol}, see~\cite{Fleming-Soner}.


\section{Regularity for the Stationary Problem.}
\label{regsection}

In order to state and prove regularity results for the stationary problem, we have to impose additional restrictions on the possible singularities of the measures $\nu_x$.

\medskip
\noindent
{\bf (M1)} \textsl{ There exists a constant $C_1 > 0$ and $\sigma \in (0,1)$ such that, for all $\beta \in [0,2]$ and all $\delta > 0$, we have
\begin{equation*}
\sup \limits_{x \in \bar{\Omega}} \int_{B_\delta^c} \min \{1, |z|^\beta \} \nu_x(dz) 
\leq C_1 h_{\beta, \sigma}(\delta),
\end{equation*}
where $h_{\beta, \sigma}(\delta)$ is defined for $\delta > 0$ as
\begin{equation*}\label{defh}
h_{\beta, \sigma}(\delta) = \left \{ \begin{array}{cl} \delta^{\beta - \sigma} \quad & \mbox{if} \ \beta < \sigma \\
|\ln(\delta)| + 1 \quad & \mbox{if} \ \beta = \sigma \\
1 \quad & \mbox{if} \ \beta > \sigma, \end{array} \right . 
\end{equation*}
and where we use the convention $|z|^\beta = 1, z \in \R^N$ when $\beta = 0$.
}

\medskip
\noindent
{\bf (M2)} \textsl{ There exists $ C_2 > 0$ and $\sigma \in (0,1)$ such that, for all $\beta \in (\sigma, 2]$ and
all $\delta \in (0,1)$ we have
\begin{equation*}
\sup \limits_{x \in \bar{\Omega}} \int_{B_\delta} |z|^\beta \nu_x(dz) \leq C_2 \delta^{\beta - \sigma}. 
\end{equation*}
}

\medskip

Roughly speaking, assumptions (M1) and (M2) say $\I$ has at most order $\sigma \in (0,1)$. Below we always assume that (M1) and (M2) are satisfied with the same $\sigma$.

Examples of censored nonlocal operator satisfying (M0)-(M2) are the censored fractional Laplacian of order $\sigma$ 
(see~\eqref{censoredLaplacian}) and regional operators depending on the distance to the boundary of the domain 
(see~\cite{Ishii-Nakamura}), typically
\begin{equation}\label{regional}
\I(u,x) = \int_{B_{d(x)}} [u(x + z) - u(x)] |z|^{-(N + \sigma)} dz.
\end{equation}

In what follows, we are going to argue on the simpler equation
\begin{equation}\label{Eq}
-\I(u,x) + b_0 |Du|^m = A_0, \quad \mbox{in}  \ \Omega,
\end{equation}
where $b_0 > 0$ and $A_0 \geq 0$. 

We start with the following preliminary result concerning the interior regularity for bounded subsolutions of~\eqref{Eq}.
\begin{lema}\label{liplocal} 
Let $\I$ be a nonlocal operator as in~\eqref{operator} satisfying (M0)-(M2) with 
the same $\sigma \in (0,1)$ and $m > \sigma$. If $u$ is a bounded viscosity subsolution to Equation~\eqref{Eq}, 
then there exists $K, r_0 > 0$ such that for all $0 < r < r_0$, we have
\begin{eqnarray*}
|u(x) - u(y)| \leq K r^{-1} |x - y|, \quad  \hbox{for all}\ x,y \in \Omega \setminus \Omega_r. 
\end{eqnarray*}

The constant $K$ only depends on the data and $\mathrm{osc}_{\bar{\Omega}}(u)$.
\end{lema}

\noindent
{\bf \textit{Proof:}} For $x \in \Omega$, we consider the function
\begin{equation*}\label{testinglemaholdersubsol}
\Phi_x: y \mapsto u(y) - \phi_x(y), \quad y \in \bar{\Omega} 
\end{equation*}
with $\phi_x(y) := u(x) + Kd^{-1}(x)|x - y|$ and $K > 0$ to be fixed later. Our aim is to prove that $\Phi_x(y) \leq 0$ for any $y \in \bar{\Omega}$ for some $K$ large enough which does not depend on $x$. Indeed, if this is true for any $x \in \Omega \setminus \Omega_r$, the result is proved. 

Note that the function $\Phi_x$ is u.s.c. on $\bar{\Omega}$ and therefore it attains its (nonnegative) maximum  
at a point $\bar{y} \in \bar{\Omega}$. We argue by contradiction assuming that $\Phi_x(\bar{y}) >0$. Since $\Phi_x (x)=0$, we clearly have $\bar{y}\neq x$.

Furthermore, if we take $K \geq 4\mathrm{osc}_{\bar{\Omega}}(u)$, 
from the inequality $\Phi_x(x) < \Phi_x(\bar{y})$ we see that $\bar{y} \in \Omega$. Indeed, it implies
$$ Kd(x)^{-1}|x - \bar{y}| \leq 2\mathrm{osc}_{\bar{\Omega}}(u)\; ,$$
and therefore $|x - \bar{y}| \leq d(x)/2$, which leads to $d(\bar{y}) \geq d(x)/2$. 

We can write the viscosity subsolution inequality for $u$ at $\bar{y}$ since $\phi_x$ is smooth enough at $\bar{y}$ and, 
for each $0 < \delta \leq d(x)/2$, we have
\begin{equation}\label{testinglemma1}
- \mathcal{I}[B_\delta](\phi_x, \bar{y}) - \mathcal{I}[B_\delta^c](u, \bar{y}) + b_0 K^m d^{-m}(x) \leq A_0.
\end{equation}

By the Lipschitz continuity of $\phi_x$ and by (M2) we readily get
\begin{equation*}
\I[B_\delta](\phi_x, \bar{y}) \leq C_2 Kd^{-1}(x)\delta^{1 - \sigma}, 
\end{equation*}
where the constant $C_2$ appears in (M2). On the other hand, by (M1) we see that
\begin{equation*}
\I[B_\delta^c](u, \bar{y}) \leq 2C_1 \mathrm{osc}_{\bar{\Omega}}(u) \delta^{-\sigma},
\end{equation*}
and then, replacing these estimates into~\eqref{testinglemma1} we obtain
\begin{equation}\label{ineqlemma1}
b_0 K^m d(x)^{-m} \leq C_2 K d(x)^{-1} \delta^{1 -\sigma}  + C_1 \mathrm{osc}_{\bar{\Omega}}(u) \delta^{-\sigma} + A_0. 
\end{equation}
and  we end up with
$$
K^m \leq C \Big{(} K d(x)^{m -1}\delta^{1-\sigma} + \mathrm{osc}_{\bar{\Omega}}(u) d(x)^{m}\delta^{-\sigma} + d(x)^{m} \Big{)},
$$
for some constant $C > 0$ depending only on the data.

Now we choose $\delta = (1+\mathrm{osc}_{\bar{\Omega}}(u))d(x)/K$ (assuming that $K>2(1+\mathrm{osc}_{\bar{\Omega}}(u))$) in the above expression which becomes
$$
K^m \leq C \Big{(} K^\sigma (1+\mathrm{osc}_{\bar{\Omega}}(u))^{1-\sigma} d(x)^{m-\sigma} + 1\Big{)},
$$

From this last inequality and since $m > \sigma$, we conclude the result by choosing $K$ large enough to get a contradiction. 
We point out that the above analysis drives us to an estimate of the size of $K$ 
like 
\begin{equation*}
C \max \{ (1+\mathrm{osc}_{\bar{\Omega}}(u))^{(1-\sigma)/(m-\sigma)}, 1+\mathrm{osc}_{\bar{\Omega}}(u) \},
\end{equation*}
for some $C > 0$ not depending on $r$ nor $\mathrm{osc}_{\bar{\Omega}}(u)$.
The proof is complete.
\qed

\medskip

The next lemma improves the Lipschitz seminorm.
\begin{lema}\label{improliplocal} 
Assume that $\Omega$ is a bounded $C^1$-domain, 
let $\I$ be a nonlocal operator as in~\eqref{operator} satisfying (M0)-(M2) with the same $\sigma \in (0,1)$ 
and $m > \sigma$. If $u$ is a bounded viscosity subsolution to Equation~\eqref{Eq}, then there exists $C, r_0 > 0$ such that
for all $0 < r < r_0$, we have
\begin{equation*}
|u(x) - u(y)| \leq \bar{C} r^{-\sigma/m} |x - y| \quad \mbox{for} \ x, y \in \Omega \setminus \Omega_r,
\end{equation*}
where the constant $\bar{C}$ depends on the data and the constant $K$ of Lemma~\ref{liplocal}.
\end{lema}

\noindent
{\bf \textit{Proof:}} From Lemma~\ref{liplocal} and Rademacher's Theorem, we know that any bounded subsolution of \eqref{Eq} is differentiable a.e. in $\Omega$.
Hence, we can evaluate the equation a.e. and therefore we can write
\begin{eqnarray*}
b_0 |Du(x)|^m \leq A_0 + \I(u,x)\quad \hbox{a.e. in  }\Omega.
\end{eqnarray*}

Writing
\begin{equation*}
\I(u,x) = \I[B_{d(x)/2}^c](u, x) + \I[B_{d(x)/2}](u, x), 
\end{equation*}
by (M1) and (M2) we arrive at
\begin{equation*}
b_0 |Du(x)|^m \leq  A_0 + C \mathrm{osc}_{\bar{\Omega}}(u) d(x)^{-\sigma} + C d(x)^{1 - \sigma} \underset{y \in B_{d(x)/2}(x)}{\mbox{essup}} \{ |Du(y)| \} ,
\end{equation*}
for some $C > 0$ not depending on $x, K$ or $\mathrm{osc}_{\bar{\Omega}}(u)$. Applying Lemma~\ref{liplocal} over $|Du(y)|$ we conclude 
\begin{eqnarray}\label{newK}
|Du(x)| \leq \bar{C} d(x)^{-\sigma/m},
\end{eqnarray}
where 
\begin{equation}\label{barC}
\bar{C} = C (1 + K)^{1/m},
\end{equation}
with $C > 0$ depending on the data, and $K$ is given by Lemma~\ref{liplocal}. 

The result is a direct consequence of this inequality. In fact,
consider $0 < r < \min \{ \delta_0, r_0\}$, where $r_0$ is given in Lemma~\ref{liplocal} and $\delta_0$ is such that the distance function is smooth in $\Omega_{\delta_0}$. We recall that, for any $x \in \Omega_{\delta_0}$, there exists a unique ``projection'' $\hat{x} \in \partial \Omega$ such that $d(x)=|x-\hat{x}|$ and $\hat{x}$ is equal to $=x-d(x)Dd(x)$. The ``$\hat{\ }$'' notation always denotes below such a projection.

Let $x, y \in \Omega \setminus \Omega_r$. If 
$|x - y| \leq r/2$ then 
$$
[x,y] := \{ t x + (1 - t)y : t \in (0,1)\} \subset \Omega.
$$ 

Moreover, for each $z \in [x,y]$ we have $d(z) \geq r/2$. In fact
\begin{equation*}
\begin{split}
d(z) = & d(tx + (1 - t)y) = d(x) + \left(d(tx + (1 - t)y)-d(x)\right)\\
\geq & d(x) - |x - (tx + (1 - t) y)| \\
\geq & d(x) - (1 - t)|x - y| \\
\geq & r - r/2(1 - t),
\end{split}
\end{equation*}
concluding $d(z) \geq r/2$ since $t \in [0,1]$. We use this in the following formal computation (arguing as if $u$ were $C^1$, but the justification is more than classical) together with~\eqref{newK} to write down
\begin{equation*}
\begin{split}
u(x) - u(y) = & \int \limits_{0}^{1} \frac{d}{dt} u(tx + (1 - t)y)dt \\
\leq & C \int \limits_{0}^{1} d(tx + (1 - t)y)^{-\sigma/m}|x - y| dt \\
\leq & C \min \{d(x), d(y)\}^{-\sigma/m} |x - y|,
\end{split}
\end{equation*}
implying that 
\begin{equation}\label{improlipsublocalI}
|u(x) - u(y)| \leq Cr^{-\sigma/m} |x - y|.
\end{equation}

On the other hand, if $|x - y| > r/2$, the only difficulty is close to the boundary and we may assume $d(x), d(y) < \delta_0$. Then, we denote $x_0, y_0 \in \Omega$ as the points such that
$\hat{x}_0 = \hat{x}, \hat{y}_0 = \hat{y}$ and $d(x_0) = d(y_0) = \delta_0$. Then, we have
\begin{equation}\label{improlipsublocalII}
u(x) - u(y) = (u(x) - u(x_0)) + (u(x_0) - u(y_0)) + (u(y_0) - u(y)). 
\end{equation}

For the first and third term on the right hand side 
of~\eqref{improlipsublocalII} we clearly see that $[x,x_0], [y,y_0] \subset \Omega$ and then we can apply the same argument 
leading to~\eqref{improlipsublocalI} to conclude
\begin{equation*}
|u(x) - u(x_0)|, |u(y) - u(y_0)| \leq Cr^{-\sigma/m}|x - y|.
\end{equation*}

On the other hand, for the second term of the right-hand side of~\eqref{improlipsublocalII} we apply Lemma~\ref{liplocal} and the regularity 
of the boundary to conclude 
$$
|u(x_0) - u(y_0)| \leq C|x - y|
$$ 
with a constant $C$ independent of $r$. This concludes the proof.
\qed


With this last lemma we are in position to prove the H\"older regularity up to the boundary.
\begin{prop}\label{propholder}\textsc{(Global H\"older Estimates)}
Assume that $\Omega$ is a bounded $C^1$-domain, let $\I$ be a nonlocal operator as in~\eqref{operator} 
satisfying (M0)-(M2) with the same $\sigma \in (0,1)$ and $m > \sigma$. If $u$ is a bounded viscosity 
subsolution to Equation~\eqref{Eq}, then there exists a constant $C_0 > 0$ such that
\begin{eqnarray*}
|u(x) - u(y)| \leq C_0|x - y|^{\frac{m - \sigma}{m}}, \quad \hbox{for all}\ x, y \in \Omega,
\end{eqnarray*}
where $C_0$ depends on the data and $\mathrm{osc}_{\bar{\Omega}}(u)$.

In particular, $u: \Omega \to \R$ can be extended as a continuous function on $\bar{\Omega}$.
\end{prop}

\noindent
{\bf \textit{Proof:}} If $|x - y| \leq \min\{ d(x), d(y)\}$, then we use directly Lemma~\ref{improliplocal} to conclude
\begin{equation*}
|u(x) - u(y)| \leq \bar{C} \min\{ d(x), d(y)\}^{- \sigma/m} |x - y| \leq C \bar{C} |x - y|^{1 - \sigma/m},
\end{equation*}
where $\bar{C}$ is given by~\eqref{barC} and $C > 0$ depends only on the data.

Now, assuming $\min \{ d(x) , d(y)\} < |x - y|$ (which implies that both $x$ and $y$ are close to the boundary since $|x-y|$ can be chosen less than (say) $\delta_0/2$), we proceed in a similar way as Lemma~\ref{improliplocal}, considering a parameter $\delta$ with
$|x - y| \leq \delta \leq \delta_0$ and define $x_\delta, y_\delta \in \Omega$ such that $d(x_\delta) = d(y_\delta) = \delta$ and 
$\hat{x}_\delta = \hat{x}, \hat{y}_\delta = \hat{y}$. Then, following~\eqref{improlipsublocalII} we write
\begin{equation*}
u(x) - u(y) = u(x) - u(x_\delta) + u(x_\delta)- u(y_\delta) + u(y_\delta) - u(y),
\end{equation*}
and using Lemma~\ref{improliplocal} we get
$
|u(x_\delta)- u(y_\delta)| \leq \bar{C} \delta^{-\sigma/m} |x - y|,
$
meanwhile, using again Lemma~\ref{improliplocal} we can write (again formally but this is easy to justify)
\begin{equation*}
\begin{split}
u(x) - u(x_\delta) \leq& \int \limits_{0}^{1} |Du(t x + (1 - t)x_\delta)| |x - x_\delta|dt\\
\leq & \bar{C} \int \limits_{0}^{1} (td(x) + (1 - t)\delta)^{-\sigma/m} (\delta - d(x)) dt\\
\leq & \bar{C} (1 - \sigma/m)^{-1} (\delta^{1-\sigma/m} - d(x)^{1-\sigma/m}),
\end{split}
\end{equation*}
and in the same way for $u(y_\delta) - u(y)$. These estimates imply the existence of a constant $C > 0$ such that
\begin{equation*}
|u(x) - u(y)| \leq C \bar{C} (\delta^{1 - \sigma/m} + \delta^{-\sigma/m} |x - y|) \quad \mbox{for all} \ \delta \geq |x - y|.
\end{equation*}

Taking infimum over these $\delta$ we arrive to the inequality
$$
|u(x) - u(y)| \leq C \bar{C}|x-y|^{\frac{m-\sigma}{m}}
$$
where $C > 0$ depends only on the data. We finish recalling that by the definition 
$\bar{C} = C (1 + K)^{1/m},$ where the constraints on $K$ are given at the end of the proof of Lemma~\ref{liplocal}.
\qed

\medskip

The first consequence of the last proposition is the following
\begin{cor}\label{oscbound}\textsc{(Oscillation Bound)}
Assume that $\Omega$ is a bounded $C^1$-domain, let $\I$ be a nonlocal operator as in~\eqref{operator} 
satisfying (M0)-(M2) with the same $\sigma \in (0,1)$ and $m > 1$. Let $u$ be a bounded viscosity subsolution 
to Equation~\eqref{Eq} and consider the function
\begin{equation*}
\tilde{u}(x) = \left \{ \begin{array}{ll} u(x) \quad & \mbox{if} \ x \in \Omega \\ 
\limsup \limits_{y \to x, y \in \Omega} u(y) \quad & \mbox{if} \ x \in \partial. \Omega \end{array} \right . 
\end{equation*}

Then, there exists $C > 0$  depending only on the data such that
\begin{equation*}
\mathrm{osc}_{\bar{\Omega}}(\tilde{u}) \leq C.
\end{equation*}
\end{cor}

\noindent
{\bf \textit{Proof:}} Note $\tilde{u} = u$ in $\Omega$. In view of Proposition~\ref{propholder}, 
we see that $\tilde{u}$ is well-defined in $\bar{\Omega}$ and it is H\"older continuous on $\bar{\Omega}$. Moreover, we can write
\begin{equation*}
|\tilde{u}(x) - \tilde{u}(y)| \leq C (1 + K)^{1/m} |x - y|^{(m - \sigma)/m},
\end{equation*}
for each $x, y \in \bar{\Omega}$. Thus, we can take $x_{min}, x_{max} \in \bar{\Omega}$ making 
\begin{equation*}
\tilde{u}(x_{max}) - \tilde{u}(x_{min}) = \mathrm{osc}_{\bar{\Omega}}(\tilde{u}),
\end{equation*}
from which we obtain
\begin{equation*}
\mathrm{osc}_{\bar{\Omega}}(\tilde{u}) \leq C (1 + K)^{1/m} \mathrm{diam}(\Omega)^{(m - \sigma)/m}. 
\end{equation*}

But the end of the proof of Lemma~\ref{liplocal} provides the constraint
$$
K = C \max \{ (1+\mathrm{osc}_{\bar{\Omega}}(u))^{(1-\sigma)/(m-\sigma)}, 1+\mathrm{osc}_{\bar{\Omega}}(u) \}
$$
for some constant $C>0$ depending only on the data. 
Since we assume $m>1$, the exponent $(1-\sigma)/(m-\sigma)$ is less than $1$ and 
therefore we can choose $K$ as $C'(1+\mathrm{osc}_{\bar{\Omega}}(u))$ for some $C'$ large enough depending only on the data.
Hence, changing perhaps the constant, the above inequality reads
$$\mathrm{osc}_{\bar{\Omega}}(\tilde{u}) \leq C (1 + \mathrm{osc}_{\bar{\Omega}}(\tilde{u}))^{1/m} \mathrm{diam}(\Omega)^{(m - \sigma)/m}, $$
providing an estimate on the oscillation of $u$.
\qed


\medskip

The following is the main result of this paper
\begin{teo}\label{teolip}\textsc{(Global Lipschitz Regularity)}
Assume that $\Omega$ is a bounded $C^1$-domain, let $\I$ be a nonlocal operator as in~\eqref{operator} 
satisfying (M0)-(M2) with the same $\sigma \in (0,1)$ and $m > \sigma$. 
If $u$ is a bounded viscosity subsolution to Equation~\eqref{Eq}, then there exists a constant $L > 0$ such that 
\begin{eqnarray*}
|u(x) - u(y)| \leq L|x - y|, \quad \mbox{for all} \ x, y \in \Omega,
\end{eqnarray*}
where the constant $L$ depends on the data and $\mathrm{osc}_{\bar{\Omega}}(u)$.
\end{teo}

We point out that, in this result, $m$ can be less or equal to $1$. 
Notice that, while for $m>1$ Corollary~\ref{oscbound} provides a bound on the oscillation of $u$, in general 
we do not have such bound when $\sigma < m \leq 1$. 
Therefore, in Theorem~\ref{teolip}, the oscillation can be seen as an external and additional data. 

In proving this theorem we require the following
\begin{lema}\label{lemaId}
Let $\I$ be a nonlocal operator as in~\eqref{operator} satisfying (M0)-(M2) with the same $\sigma \in (0,1)$. 
Then, there exists $\bar{\delta} \in (0, \delta_0)$ and $C > 0$ such that,
for each $\beta \in (0, \sigma)$
\begin{equation*}
\I(d^\beta, x) \leq Cd^{\beta - \sigma}(x), \quad \mbox{for all} \ x \in \Omega_{\bar{\delta}}. 
\end{equation*}
\end{lema}

\noindent
{\bf \textit{Proof:}} Note that $d^\beta: \bar{\Omega} \to \R$ is bounded and smooth in $\Omega_{\delta_0}$. We start considering $\bar{\delta} < \delta_0/4$
and for $x \in \Omega_{\bar{\delta}}$ we write
\begin{equation}\label{lemaId1}
\begin{split} 
\I(d^\beta, x) = \I[B_{d(x)/2}](d^\beta, x) + \I[B_{\delta_0/4} \setminus B_{d(x)/2}](d^\beta, x) + \I[B_{\delta_0/4}^c](d^\beta, x),
\end{split}
\end{equation}
and we estimate each term separately. For the first integral term in the right-hand side of~\eqref{lemaId1}, we can perform a Taylor expansion 
to write
\begin{equation*}
\I[B_{d(x)/2}](d^\beta, x) \leq \beta (d(x)/2)^{\beta - 1} \int_{B_{d(x)/2}} |z| \nu_x(dz),
\end{equation*}
from which, applying (M2) we arrive at
\begin{equation*}
\I[B_{d(x)/2}](d^\beta, x) \leq C d^{\beta - \sigma}(x).
\end{equation*}

For the second integral term in~\eqref{lemaId1}, we use that $z \mapsto d^\beta(z)$ is $\beta$-H\"older continuous in $\Omega_{\delta_0}$ and therefore
we can write
\begin{equation*}
\I[B_{\delta_0/4} \setminus B_{d(x)/2}](d^\beta, x) \leq C \int_{B_{\delta_0/4} \setminus B_{d(x)/2}} |z|^\beta \nu_x(dz),
\end{equation*}
with $C > 0$ not depending on $x$ or $\beta$. Applying (M1), we conclude
\begin{equation*}
\I[B_{\delta_0/4} \setminus B_{d(x)/2}](d^\beta, x) \leq C d^{\beta - \sigma}(x). 
\end{equation*}

Finally, using the boundedness of $d^\beta$ in $\Omega$, (M0) and (M1), for the third integral term in the right-hand side of~\eqref{lemaId1}, 
we conclude
\begin{equation*}
\I[B_{\delta_0/4}^c](d^\beta, x) \leq C \delta_0^{-\sigma}, 
\end{equation*}
where $C > 0$ depends only on $\mathrm{diam}(\Omega)$. Gathering the above estimates in the right-hand side of~\eqref{lemaId1}, we conclude the result
taking $\bar{\delta}$ smaller if it is necessary.
\qed

\medskip
\noindent
{\bf \textit{Proof of Theorem~\ref{teolip}:}} Denote $\gamma = 1 - \sigma/m$ the H\"older exponent of Lemma~\ref{propholder}.
Let $x \in \Omega$ and $\beta \in (0, \min \{ \gamma, \sigma \})$ fixed. For $L > 1$ and $0<\eta \ll 1$, we consider the function
\begin{equation*}
y \mapsto \Phi_x(y) := u(y) - \phi_x(y),
\end{equation*}
where the function $\phi_x$ has the form
\begin{equation*}
\phi_x(y) = L |y - x| - \eta d^\beta(y).
\end{equation*}

Our aim is to prove that, for $L$ large enough (with a size which does not depend on $x$), we have $\Phi_x (y)\leq 0$ on $\bar{\Omega}$, for all $\eta$ small enough (possibly depending on $L$ and/or $x$). If this is true, we deduce the Lipschitz continuity from this property by letting $\eta$ tend to $0$. Of course, the difficulty is to get such property close to the boundary and therefore we may assume that $d(x) \leq \bar{\delta}/2$ with $\bar{\delta}$ as in Lemma~\ref{lemaId}.

We argue by contradiction, assuming that
\begin{equation}\label{maxlip}
\max \limits_{\bar{\Omega}} \Phi_x \geq \epsilon_L > 0,
\end{equation}
for each $L$ large and $\eta$ small.

By the upper semicontinuity of $u$, there exists $\bar{y} \in \bar{\Omega}$ realizing the maximum in~\eqref{maxlip}, and this  maximum cannot be equal to $x$ since $\epsilon_L > 0$. Moreover, $\bar{y} \in \Omega$ since otherwise, for $a > 0$ small enough, the inequality $\Phi(\bar{y} + aDd(\bar{x})) \leq \Phi_x(\bar{y})$ leads to
\begin{equation*}
\eta a^\beta \leq u(\bar{y}) - u(\bar{y} + aDd(\bar{y})) + L (|\bar{y} + aDd(\bar{y}) - x| - |\bar{y} - x|),
\end{equation*}
and by the H\"older regularity given by Proposition~\ref{propholder}, we arrive at
\begin{equation*}
\eta a^\beta \leq C_0 a^\gamma + La, 
\end{equation*}
which is a contradiction when $a > 0$ is small since $\beta <\gamma$.

Since $\bar{y}\in \Omega$ and $\bar{y} \neq x$, we can write the viscosity subsolution inequality for $u$ with 
test function $\phi_x$ at $\bar{y}$. Thus, for all $\delta > 0$ we can write
\begin{equation}\label{testteolip}
- \I[B_\delta](\phi_x, \bar{y}) - \I[B_\delta^c](u, \bar{y}) + b_0 |D\phi_x(\bar{y})|^m \leq A_0,
\end{equation}
and the idea is to estimate each term in the left-hand side separately. Note that
\begin{equation*}
\I[B_\delta](\phi_x, \bar{y}) = LJ_1 - \eta J_2,
\end{equation*}
with
\begin{equation*}
\begin{split}
J_1 = & \int_{B_\delta} (|\bar{y} + z - x| - |\bar{y} - x|) \nu_{\bar{y}}(dz), \\
J_2 = & \int_{B_\delta} [d^\beta(\bar{y} + z) - d^\beta(\bar{y})] \nu_{\bar{y}}(dz).
\end{split}
\end{equation*}

Since $\sigma < 1$ and applying (M1), we easily conclude that 
\begin{equation*}
J_1 \leq C_1 \delta^{1 - \sigma}.
\end{equation*}

Now, using that $\Phi_x(x) \leq \Phi_x(\bar{y})$ and considering $\eta$ small in terms of $\mathrm{diam}(\Omega)$, we see that
\begin{equation}\label{x-ylip}
|\bar{y} - x| \leq L^{-1} (\mathrm{osc}_{\bar{\Omega}} (u)+ 1). 
\end{equation}

Since we are considering $x$ such that $d(x) \leq \bar{\delta}/2$ with $\bar{\delta}$ as in Lemma~\ref{lemaId}, 
by the last inequality
we can take $L$ large enough (depending on $\mathrm{osc}_{\bar{\Omega}}(u)$ and $\bar{\delta}$, and, a fortiori, on $\delta_0$) 
to get $d(\bar{y}) \leq \bar{\delta}$. Then, by Lemma~\ref{lemaId}, we can write
\begin{equation*}
J_2 \leq C d^{\beta - \sigma}(\bar{y}), 
\end{equation*}
where $C$ does not depend on $\bar{x}, y, L$ or $\eta$. Then, by the estimates for $J_1$ and $J_2$ we get
\begin{equation*}
\I[B_\delta](\phi_x, \bar{y}) \leq C L \delta^{1 - \sigma} + C \eta d^{\beta - \sigma}(\bar{y}).
\end{equation*}

Concerning the integral term outside $B_\delta$, the H\"older continuity of $u$ on $\bar{\Omega}$ given by Lemma~\ref{propholder} allows us 
to write
\begin{equation*}
\I[B_\delta^c](u, \bar{y}) \leq C \int_{B_\delta^c} |z|^\gamma \nu_{\bar{y}}(dz) 
\leq C \delta^{\gamma - \sigma},
\end{equation*}
where $C > 0$ depends only on the data and the constant $C_0$ of Proposition~\ref{propholder}. In summary, 
concerning the integral terms we obtain the estimate
\begin{equation}\label{nonlocallip}
\I[B_\delta](\phi_x, \bar{y}) + \I[B_\delta^c](u, \bar{y}) 
\leq C\Big{(} L \delta^{1 - \sigma} + \eta d^{\beta - \sigma}(\bar{x}) + \delta^{\gamma - \sigma} \Big{)}. 
\end{equation}

Now we deal with the first-order term. A straightforward computation gives us
\begin{equation*}
D\phi_x(\bar{y}) = L  \widehat{(\bar{y} - x)}- \eta \beta d^{\beta - 1}(\bar{y}) Dd(\bar{y}), 
\end{equation*}
with 
$$
\widehat{(\bar{y} - x)}:=\frac{(\bar{y} - x)}{|\bar{y} - x|}\; .
$$

Recalling that $\bar{y}$ depends on $\eta$, at this point we have two different situations : either there exists a 
subsequence such that $d(\bar{y}) \geq d(x)$ as $\eta$ tends to zero, or $d(\bar{y}) \leq d(x)$ as $\eta$ tends to zero.
In the former case, $\eta d(\bar{y})^{\beta - 1} \leq \eta d(x)^{\beta - 1}$ 
and therefore $\eta d(\bar{y})^{\beta - 1}\to 0$ as $\eta \to 0$. This yields $|D\phi_x(\bar{y})| \geq L +o_\eta(1)$.
In the second case, computing  $|D\phi_x(\bar{y})|^2$, we see that
\begin{equation*}
|D\phi_x(\bar{y})|^2 =  L^2 + \eta^2 \beta^2 d^{2(\beta - 1)}(\bar{y}) 
- 2 \eta \beta L d^{\beta - 1}(\bar{y}) \langle Dd(\bar{y}), \widehat{(\bar{y} - x)} \rangle.
\end{equation*}

But $x,\bar{y} \in \Omega_{\delta_0}$ where $d$ is $C^1$, and performing a Taylor expansion we get
\begin{equation*}
d(x) - d(\bar{y}) =  \langle Dd(\bar{y}) , (x - \bar{y}) \rangle + o(|\bar{y} - x|),
\end{equation*}
which replaced into the above expression for $|D\phi_x(\bar{x})|^2$ drives us to
\begin{equation}\label{gradientlip0}
\begin{split}
|D\phi_x(\bar{y})|^2 = & L^2 + \eta^2 \beta^2 d^{2(\beta - 1)}(\bar{y}) \\
& - 2 \eta \beta L d^{\beta - 1}(\bar{y}) ((d(\bar{y}) - d(x))/|\bar{y} - x| - o(1)),
\end{split}
\end{equation}
where $o(1)\to 0$ when $|\bar{x} - y|\to 0$. But taking into account that $d(\bar{y}) \leq d(x)$, we are led to
$$
|D\phi_x(\bar{y})|^2 \geq  L^2 + \eta^2 \beta^2 d^{2(\beta - 1)}(\bar{y})
- 2 \eta \beta L d^{\beta - 1}(\bar{y}) o(1),
$$
and by using Cauchy-Schwarz inequality, we arrive at
\begin{equation}\label{gradientlip}
|D\phi_x(\bar{y})|^2 \geq  L^2/2 + \eta^2 \beta^2 d^{2(\beta - 1)}(\bar{y})/2,
\end{equation}
for all $L$ large and $\eta$ small enough.

Then, depending on the case we are looking at, we replace~\eqref{nonlocallip} and~\eqref{gradientlip} or $|D\phi_x(\bar{y})| \geq L +o_\eta(1)$ into~\eqref{testteolip} to obtain, for all $x \in \Omega, \delta > 0$, 
$L$ large (depending only on the data and the oscillation of $u$) and $\eta$ small depending on $L$ and $x$, the inequality
\begin{equation}\label{ineqlip}
L^m + (\eta d^{\beta- 1}(\bar{y}))^m 
\leq C \Big{(} A_0 + L \delta^{1 - \sigma} + \delta^{\gamma - \sigma} + \eta d^{\beta - \sigma}(\bar{y}) \Big{)},
\end{equation}
for some constant $C > 0$ depending only on the data. 

From this inequality, we claim that there exists an universal constant $\Lambda_0 > 0$ such that
\begin{equation}\label{claimlip}
C\eta d^{\beta - \sigma}(\bar{y}) - (\eta d^{\beta- 1}(\bar{y}))^m \leq \Lambda_0,
\end{equation}
for all $\eta$ small. We postpone the justification of this claim until the end of this proof. 

Using this claim into~\eqref{ineqlip}, we observe that
\begin{equation*}
L^m \leq C \Big{(} A_0 + L \delta^{1 - \sigma} + \delta^{\gamma - \sigma} + \Lambda_0 \Big{)}.
\end{equation*}

Hence, in the case $m \geq 1$, since $\sigma < 1$ we can fix $\delta > 0$ small in terms of $C$ to conclude the result by taking $L$ large in terms 
of the data and the oscillation of $u$. In the case $m < 1$, we take $\delta = L^{-1} > 0$ and since $m > \sigma$, we conclude taking 
$L$ large in terms of the data and the oscillation of $u$. This concludes the proof of the theorem, up to the justification 
of~\eqref{claimlip} which is provided next.

We write
\begin{equation*}
\theta = \eta d^{\beta - \sigma}(\bar{y}), 
\end{equation*}
from which, we obtain that
\begin{equation*}
C\eta d^{\beta - \sigma}(\bar{y}) - (\eta d^{\beta- 1}(\bar{y}))^m 
= C \theta - \theta^{\tau} \eta^{-\alpha},
\end{equation*}
where $\tau := m(1 - \beta)/(\sigma - \beta) > 0$ and $\alpha := m(1 - \sigma)/(\sigma - \beta) > 0$. 
Note that since $m > \sigma$, we can fix $\beta > 0$ small enough to get $\tau > 1$. With this choice we can consider
$\eta \leq C^{-1/\alpha}$ and with this we arrive at
\begin{equation*}
C\eta d^{\beta - \sigma}(\bar{y}) - (\eta d^{\beta- 1}(\bar{y}))^m  \leq C (\theta - \theta^{\tau}) 
\leq C \tau^{-\tau/(\tau - 1)} (\tau - 1) \leq C (\tau - 1) =: \Lambda_0,
\end{equation*}
which is a constant depending only on the data, but not on $\eta, L$ or $x$. This concludes the claim.
\qed


\section{Comparison Principle and Well-Posedness.}\label{cpw}

Here we address the well-posedness of the following Cauchy-Dirichlet evolution problem
\begin{equation}\label{eqparabolic}\tag{CP}
\left \{ \begin{array}{rll} u_t - \I(u(\cdot, t), x) + H(x, u, Du) &= f, \ &\mbox{in} \ Q \\
u &= \varphi, \ &\mbox{in} \ \partial Q \\
u &= u_0, \ &\mbox{in} \ \bar{\Omega} \times \{ 0 \}. \end{array} \right .
\end{equation}
where $\I$ is a nonlocal operator with the form~\eqref{operator} satisfying conditions (M0)-(M2),
$H \in C(\bar{\Omega} \times \R \times \R^N)$, $f \in C(\bar{Q})$, $\varphi \in C(\partial Q)$ 
and $u_0 \in C(\bar{\Omega})$. We always assume that \eqref{compatibility} holds.

On $H$ we additionally assume the properness condition

\medskip
\noindent
{\bf (H1)} \textsl{For all $R > 0$, there exists $\lambda_R \geq 0$ such that, for all 
$x \in \bar{\Omega}$, $u\geq v$ with $|u|, |v| \leq R$ and $p \in \R^N$ we have}
\begin{equation*}
H(x,u,p) - H(x,v,p) \geq \lambda_R  (u - v).
\end{equation*}

In order to apply  the regularity results of the previous section, we consider the coercivity condition

\medskip
\noindent
{\bf (H2)} \textsl{There exists $m > 0$ and $\underline{c} > 0$ such that, for each $R > 0$
there exists $C_R \geq 0$ satisfying
\begin{equation*}
H(x, r, p) \geq \underline{c}|p|^m - C_R,
\end{equation*}
for all $x \in \bar{\Omega}, p \in \R^N$ and $|r| \leq R$.
}

\medskip

Finally, in addition to (M0)-(M2), we will require the following continuity assumption over the 
measures $\{ \nu_x \}$ defining $\I$

\medskip
\noindent
{\bf (M3)} \textsl{Assume (M1), (M2) hold with the same value $\sigma \in (0,1)$ and let $h$ be as in (M1). 
Then, for each $\beta > 0$, 
there exists a modulus of continuity $\omega_{\beta, \sigma}$ such that for all $x, y \in \bar{\Omega}$ and $\delta > 0$, we have
\begin{equation*}
\int_{B_\delta^c} |z|^\beta|\nu_x(dz) - \nu_y(dz)| \leq (1 + h_{\beta, \sigma}(\delta)) \ \omega_{\beta, \sigma}(|x - y|)
\end{equation*}
}
where $|\nu_x(dz) - \nu_y(dz)|$ denotes the total variation of the measure $\nu_x(dz) - \nu_y(dz)$.

Note that (M3) holds in the censored fractional Laplacian case (see~\eqref{censoredLaplacian})
and regional operators depending on the distance to the boundary (see~\eqref{regional}).

\medskip

The main result of this section reads as follows
\begin{prop}\label{comparison}
Assume that $\Omega$ is a bounded, $C^1$-domain and let $\I$ be a nonlocal operator as in~\eqref{operator} satisfying (M0)-(M3) with the same $\sigma \in (0,1)$. Assume also that $H\in C(\bar{\Omega} \times \R \times \R^N)$ satisfies (H1)-(H2) with $m > \sigma$, $f \in C(\bar{Q})$, $\varphi \in C(\partial Q)$, 
$u_0 \in C(\bar{\Omega})$ and that~\eqref{compatibility} holds. 
If $u, v$ are respectively viscosity subsolution and supersolution to~\eqref{eqparabolic} which are 
bounded in each compact subset of $\bar{Q}$, then
\begin{equation*}
u \leq v \quad \mbox{in} \ Q \cup \bar{\Omega} \times \{ 0 \}.
\end{equation*}

Moreover, if $\tilde{u}$ is defined by
\begin{equation}\label{tildeu}
\tilde{u}(x,t) = \left \{ \begin{array}{ll} u(x,t) \quad & \mbox{if} \ (x, t) \in Q \cup \bar{\Omega} \times \{ 0 \} \\
\limsup \limits_{(y, s) \to (x,t), y \in \Omega} u(y,s) \quad & \mbox{if} \ (x, t) \in \partial Q,
\end{array} \right .
\end{equation}
then $\tilde{u} \leq v$ in $\bar{Q}$.
\end{prop}


Before giving the proof of this result, we start with two lemmas concerning the (parabolic) boundary condition of problem~\eqref{eqparabolic}.
Next result states classical initial condition holds mainly by the compatibility condition~\eqref{compatibility}
\begin{lema}\label{lemat_0}Under the assumptions of Proposition~\ref{comparison}, for all $x \in \bar{\Omega}$ we have $u(x, 0) \leq u_0(x) \leq v(x, 0)$.
\end{lema}

We refer to~\cite{DaLio} for a proof of Lemma~\ref{lemat_0} in the second-order case which can be readily adapted to the current framework. We continue with the classical boundary condition
for subsolutions
\begin{lema}\label{lemasubboundary}
Under the assumptions of Proposition~\ref{comparison}, $u(x,t) \leq \varphi(x, t)$ for each $(x, t) \in \partial Q$. 
\end{lema}

We may refer to~\cite{Barles-Chasseigne-Imbert, BT, Topp} for a proof of this result in very similar settings.

\begin{remark}\label{rmktildeu}
A classical strategy to deal with the well-posedness of~\eqref{eqparabolic} in $C(\bar{Q})$
is to argue over the redefined function $\tilde{u}$ in~\eqref{tildeu} instead of the original subsolution $u$, 
see~\cite{Barles-Chasseigne-Imbert, BT, Topp} and the second-order references therein.
Note that $\tilde{u} = u$ in $Q$, $\tilde{u} \leq u$ on $\partial \Omega \cup \bar{\Omega} \times \{ 0 \}$ and therefore
Lemmas~\ref{lemat_0} and~\ref{lemasubboundary} hold for $\tilde{u}$. Moreover, $\tilde{u}$ is a (generalized)
viscosity subsolution to the problem if $u$ is. Thus, 
for simplicity, we are going to assume that $u = \tilde{u}$ on $\bar{Q}$ in order to avoid the superscript ``$\sim$''.
\end{remark}

The arguments to come are carried out on the finite time horizon problem
\begin{equation}\label{parabolicT}\tag{$\mathrm{CP_T}$}
\left \{ \begin{array}{rll} \partial_t u - \I(u(\cdot, t), x) + H(x, u, Du) &= f, \ &\mbox{in} \ Q_T \\
u &= \varphi, \ &\mbox{in} \ \partial Q_T \\
u &= u_0, \ &\mbox{in} \ \bar{\Omega} \times \{ 0 \}, \end{array} \right .
\end{equation}
where, for $T > 0$, we denote $Q_T = \Omega \times (0,T]$ and $\partial Q_T = \partial \Omega \times (0, T]$.
The infinite horizon setting for Proposition~\ref{comparison} is readily obtained by taking a sequence of 
problems~\eqref{parabolicT}  with $T \to +\infty$.

We require the following lemmas in order to apply the regularity results of the previous section
\begin{lema}\label{lemaugamma}
Assume $\varphi \in C(\partial Q_T), f \in C(\bar{Q}_T)$ and $H \in C(\bar{\Omega} \times \R \times \R^N)$ satisfies (H2).
Let $u$ be bounded u.s.c. viscosity subsolution to problem
\begin{eqnarray}\label{eqtolinearize}
\left \{ \begin{array}{rll} \partial_t u - \I(u, x) + H(x,u,Du) & = f \quad & \mbox{in} \ Q_T \\
u & = \varphi \quad & \mbox{in} \ \partial Q_T, \end{array} \right .
\end{eqnarray}

For $\gamma > 0$ and $(x, t) \in \bar{Q}_T$,
consider $u^\gamma$, the sup-convolution in time of $u$, which is given by the expression
\begin{equation}\label{ugamma}
u^\gamma(x,t) = \sup \limits_{s \in [0,T]} \{ u(x,s) - \gamma^{-1} (s - t)^2 \}. 
\end{equation}

Then, there exists a constant $a_\gamma > 0$ with $a_\gamma  \to 0$ as $\gamma \to 0$ such that $u^\gamma$ is a viscosity subsolution to 
problem
\begin{equation*}\label{equgamma}
\left \{ \begin{array}{rll} u^\gamma_t - \I(u^\gamma(\cdot, t), x) + H(x, u^\gamma, Du^\gamma) &= f + o_\gamma(1), \ &\mbox{in} \ \Omega \times (a_\gamma, T] \\
u^\gamma &= \varphi^\gamma, \ &\mbox{in} \ \partial \Omega \times [a_\gamma, T],  \end{array} \right . 
\end{equation*}
where $o_\gamma(1) \to 0$ as $\gamma \to 0$ and depends exclusively on the modulus of continuity of $f$.
\end{lema}

The proof of the above lemma follows closely the arguments of~\cite{Barles-Mitake, BT}. 
Note that by Lemma~\ref{lemasubboundary} the boundary condition in the above equation is satisfied in the classical sense.	

It is a well-known fact that for each $\gamma > 0$ and $x \in \bar{\Omega}$, $t \mapsto u^\gamma(x, t)$ is 
Lipschitz continuous in $[0,T]$, with Lipschitz constant $C_\gamma := 4T \gamma^{-1}$. Therefore, for any $t \in (a_\gamma, T]$, 
the function $x \mapsto u^\gamma(x, t)$ is a viscosity subsolution to
$$ - \I(u, x) + H(x, u, Du) = f + o_\gamma(1)+C_\gamma \ \mbox{in} \ \Omega \; .$$

Notice that thanks to (H2), a bounded viscosity subsolution to the above problem is a viscosity subsolution to
$$ - \I(u, x) + \underline{c} |Du|^m = f + o_\gamma(1)+ C_\gamma + C_u \ \mbox{in} \ \Omega \; ,$$
where $C_u$ depends on the $L^\infty$ norm of $u$. Thus, applying Theorem~\ref{teolip} we conclude the following
\begin{lema}\label{regularityugamma}
Let $u$ be a bounded viscosity subsolution to problem~\eqref{eqtolinearize}. 
Let $\gamma > 0$, $u^\gamma$ defined as in~\eqref{ugamma} 
and $a_\gamma$ the constant given in Lemma~\ref{lemaugamma}. Then, $u^\gamma \in \mathrm{Lip}(\bar{\Omega} \times [a_\gamma, T])$.
\end{lema}

Now we are in position to provide the

\medskip
\noindent
{\bf \textit{Proof of Proposition~\ref{comparison}:}} By contradiction, we assume that $u - v$ is positive at 
some point on $\bar{Q}_T$. Then, for all $\eta > 0$ small in terms of $T$ and the supremum of $u - v$ on $\bar{Q}_T$, we see that
\begin{equation*}\label{Mcomparison}
M := \sup \limits_{\bar{Q}_T} \{ u - v - \eta t \} > 0.
\end{equation*}

By semicontinuity, this supremum is attained at some point $(x_0, t_0) \in \bar{Q}_T$. Taking $\eta$ 
smaller if it is necessary, by Lemma~\ref{lemat_0} we have $t_0 > 0$. From this point, we fix such an $\eta > 0$.

For $\alpha > 0$, notice the function
\begin{equation*}
(x,t) \mapsto u(x,t) - v(x,t) - \eta t - \alpha(|x - x_0|^2 + (t - t_0)^2)
\end{equation*}
attains its unique maximum on $\bar{Q}_T$ at $(x_0, t_0)$, and in fact this maximum equals $M$.
Recalling the definition of $u^\gamma$ in~\eqref{ugamma}, we define 
\begin{equation}\label{Mgamma}
M_\gamma :=  \sup \limits_{(x,t) \in \bar{Q}_T} \{ u^\gamma(x,t) - v(x,t) - \eta t - \alpha(|x - x_0|^2 + (t - t_0)^2) \}.
\end{equation}

By definition of $u^\gamma$ we see that $M_\gamma \geq M$ and this supremum 
is attained at some point $(x_\gamma, t_\gamma) \in \bar{Q}_T$.
Since $(x_0, t_0)$ is a strict maximum point, properties of the sup-convolution imply that
$(x_\gamma, t_\gamma) \to (x_0, t_0)$ as $\gamma \to 0$. 

We are going to consider the case $x_\gamma \in \partial \Omega$ for all $\gamma$. The case $x_\gamma \in \Omega$
with $d(x_\gamma)$ uniformly positive with respect to $\gamma$ is obtained by easier arguments and computations, 
see~\cite{BT, Topp}.

We follow Soner's penalization procedure introduced in~\cite{Soner} (see also~\cite{usersguide, BT, Topp}).
We double variables, consider a parameter $\epsilon > 0$ and define the function
\begin{equation*}
\Phi(x,y, s, t) := u^\gamma(x, s) - v(y, t) - \phi(x,y,s,t),
\end{equation*}
where, denoting $\bar{\xi} = (Dd(x_\gamma), 0) \in \R^{N + 1}$, we define 
\begin{equation*}
\phi(x,y,s,t) := |\epsilon^{-1} ((x, s) - (y, t)) - \bar{\xi}|^2 + \eta s + \alpha(|y - x_0|^2 + (t - t_0)^2).
\end{equation*}

Since $\Phi$ is upper semicontinuous
in $\bar{Q}_T \times \bar{Q}_T$, there exists a point $(\bar{x}, \bar{y}, \bar{s}, \bar{t}) \in \bar{Q}_T \times \bar{Q}_T$ 
where the function $\Phi$ attains its maximum. By the inequality 
$$
\Phi(\bar{x}, \bar{y}, \bar{s}, \bar{t}) \geq \Phi(x_\gamma + \epsilon Dd(x_\gamma), x_\gamma, t_\gamma, t_\gamma)
$$ 
and the continuity of $u^\gamma$, classical viscosity arguments lead us to
\begin{equation}\label{comparisoneps}
\begin{split}
& (\bar{x}, \bar{s}), (\bar{y}, \bar{t}) \to (x_\gamma, t_\gamma), \quad |\epsilon^{-1}(\bar{x} - \bar{y}) - Dd(x_\gamma)| \to 0, 
\quad \epsilon^{-1}|\bar{s} - \bar{t}| \to 0 ,\\ 
& \mbox{and} \quad u^\gamma(\bar{x}, \bar{s}) \to u^\gamma(x_\gamma, t_\gamma), \quad v(\bar{y}, \bar{t}) \to v(x_\gamma, t_\gamma),
\end{split}
\end{equation}
as $\epsilon \to 0$. Moreover, for all $\epsilon$ small, the same inequality implies that
\begin{equation*}
\Phi(\bar{x}, \bar{y}, \bar{s}, \bar{t}) \geq M/2.
\end{equation*}

This inequality implies that $u^\gamma(\bar{x}, \bar{s}) > v(\bar{y}, \bar{t})$ and therefore, using
Lemma~\ref{lemasubboundary} and the continuity of $\varphi$, we deduce that, if $\bar{y} \in \partial \Omega$,
\begin{equation}\label{v<varphi}
\varphi(\bar{y}, \bar{t}) > v(\bar{y}, \bar{t}).
\end{equation}

Thus, the supersolution inequality holds at $(\bar{y}, \bar{t})$ for $v$ even if $\bar{y} \in \partial \Omega$.
On the other hand, the second convergence in~\eqref{comparisoneps} leads us to
\begin{equation}\label{xin}
\bar{x} = \bar{y} + \epsilon Dd(x_\gamma) + o(\epsilon), 
\end{equation}
form which we conclude that $\bar{x} \in \Omega$ for all $\epsilon$ small. Indeed, by a simple Taylor's expansion, $d$ being $C^1$
$$d(\bar{x}) = d(\bar{y}) + \epsilon Dd(x_\gamma) \cdot Dd(\bar{y}) + o(\epsilon),$$
and since $\bar{y}\to x_\gamma$ as $\epsilon \to 0$, we have 
$$d(\bar{x}) \geq \epsilon |Dd(x_\gamma)|^2 + o(\epsilon),$$
and $d(\bar{x})>0$ for $\epsilon$ small enough since $ |Dd(x_\gamma)|=1$.
In conclusion, $\bar{x} \in \Omega$ and we can write down a subsolution viscosity inequality for $u^\gamma$ at $(\bar{x}, \bar{s})$.

Now we substract the corresponding viscosity inequalities for $u^\gamma$ at $(\bar{x}, \bar{s})$ (see Lemma~\ref{lemaugamma}) 
and $v$ at $(\bar{y}, \bar{t})$ and then, for all $\delta > 0$ we can write
\begin{equation}\label{testinglinearization}
\mathcal{A} - \I^{\delta} - \I_{\delta} \leq o_\gamma(1),
\end{equation}
where 
\begin{equation*}
\begin{split}
\I_{\delta} = & \ \I[B_{\delta}](\phi(\cdot, \bar{y}, \bar{s}, \bar{t}), \bar{x}) 
- \I[B_{\delta}](-\phi(\bar{x}, \cdot, \bar{s}, \bar{t}), \bar{y}), \\
\I^{\delta} = & \ \I[B_{\delta}^c](u^\gamma(\cdot, \bar{s}),\bar{x}) - \I[B_{\delta}^c](v(\cdot, \bar{t}),\bar{y}), 
\end{split}
\end{equation*}
and 
\begin{equation*}
\begin{split}
\mathcal{A} = & \ (\partial_s \phi + \partial_t \phi)(\bar{x},\bar{y}, \bar{s}, \bar{t}) - f(\bar{x}, \bar{s}) + f(\bar{y}, \bar{t}) \\
& \ + H(\bar{x}, u^\gamma(\bar{x}, \bar{s}), \bar{p})- H(\bar{y}, v(\bar{y}, \bar{t}), \bar{q}),
\end{split}
\end{equation*}
with 
\begin{equation*}
\begin{split}
\bar{p} = & \ D_x \phi(\bar{x},\bar{y}, \bar{s}, \bar{t}) = \epsilon^{-1} (\epsilon^{-1}((\bar{x}, \bar{s}) - (\bar{y}, \bar{t})) 
- \bar{\xi}), \\
\bar{q} = & \ -D_y \phi(\bar{x},\bar{y}, \bar{s}, \bar{t}) = \bar{p} - 2\alpha(\bar{y} - x_0).
\end{split}
\end{equation*}

The core of the remaining proof is to estimate $\I_\delta, \I^\delta$ and $\mathcal{A}$. We start with the latterer.
Recalling that $u^\gamma \in \mathrm{Lip}(\bar{\Omega} \times [a_\gamma, T])$ and that $(\bar{x}, \bar{y}, \bar{s}, \bar{t})$
is a maximum point for $\Phi$, we have $|\bar{p}| \leq L_\gamma$ for all $\epsilon$ small, 
where $L_\gamma > 0$ is the Lipschitz constant for $u^\gamma$
given in Lemma~\ref{regularityugamma}. Denoting $R = ||u^\gamma||_\infty + ||v||_\infty$, 
applying (H1) together with the (uniform) continuity of $f$ in $\bar{Q}_T$ and~\eqref{comparisoneps}, 
we see that
\begin{equation*}
\mathcal{A} \geq \eta - o_\epsilon(1) + \lambda_R (u^\gamma(\bar{x}, \bar{s}) - v(\bar{y}, \bar{t})) 
+ H(\bar{x}, v(\bar{y}, \bar{t}), \bar{p}) - H(\bar{y}, v(\bar{y}, \bar{t}), \bar{q}),
\end{equation*}
with $o_\epsilon(1) \to 0$ as $\epsilon \to 0$ depending only on the modulus of continuity of $f$.
Hence, by the uniform continuity of $H$ on compact sets of $\bar{\Omega} \times \R \times \R^N$ we arrive at
\begin{equation*}
\mathcal{A} \geq \eta - o_\alpha(1) - o_\epsilon(1) + \lambda_R (u^\gamma(\bar{x}, \bar{s}) - v(\bar{y}, \bar{t})),
\end{equation*}
and by the definition of $M_\gamma$ in~\eqref{Mgamma} 
together with~\eqref{comparisoneps} and the continuity of $\lambda_R$, we conclude
\begin{equation}\label{A}
\mathcal{A} \geq \eta - o_\alpha(1) - o_\epsilon(1) + \lambda_R M.
\end{equation}

Now we address the nonlocal terms. By the smoothness of $\phi$ we can write
\begin{equation}\label{I_delta}
\I_{\delta} \leq \epsilon^{-1} o_{\delta}(1), 
\end{equation}
where $o_{\delta}(1) \to 0$ as $\delta \to 0$ and does not depend on $\epsilon$.

On the other hand, for $\I^{\delta}$ we can write
\begin{equation*}
\begin{split}
\I^{\delta} = & \ \int_{B_{\delta}^c} [u^\gamma(\bar{x} + z, \bar{s}) - u^\gamma(\bar{x}, \bar{s})] (\nu_{\bar{x}}(dz) - \nu_{\bar{y}}(dz)) \\
& \ + \int_{B_{\delta}^c} [u^\gamma(\bar{x} + z, \bar{s}) - v(\bar{y} + z, \bar{t}) - (u^\gamma(\bar{x}, \bar{s}) - v(\bar{y}, \bar{t}))] \nu_{\bar{y}}(dz) \\
=: & \ \I^{\delta}_1 + \I^{\delta}_2.
\end{split}
\end{equation*}

For $\I^{\delta}_1$, by Lemma~\ref{regularityugamma} and (M3) we have
\begin{equation*}
\I^{\delta}_1 \leq C \ L_\gamma \ \omega_{1, \sigma}(|\bar{x} - \bar{y}|), 
\end{equation*}
for some $C > 0$ not depending on $\epsilon$ or $\delta$. 
Thus, by~\eqref{comparisoneps} we can write
\begin{equation*}
\I^{\delta}_1 \leq L_\gamma o_\epsilon(1). 
\end{equation*}

Now, for $\I^{\delta}_2$, we divide the region of integration as 
\begin{equation*}
B_{\delta}^c = (B_{\delta}^c \cap (\Omega - \bar{x}) ) \cup (B_{\delta}^c \setminus (\Omega - \bar{x})) =: \Theta_1 \cup \Theta_2.
\end{equation*}

At this point we notice that (M0) says that $\mathrm{supp} \{ \nu_{\bar{y}}\} \subset \bar{\Omega} - \bar{y}$
and therefore, by using that $(\bar{x}, \bar{y}, \bar{s}, \bar{t})$ is a maximum for $\Phi$, we can write
\begin{equation*}
\I^{\delta}_2 
\leq C\alpha \int_{\Theta_1} |z| \nu_{\bar{y}}(dz) + 2R \int_{\Theta_2} \nu_{\bar{y}}(dz).
\end{equation*}

By (M1) we can control the first integral term in the right-hand side of the last inequality by $C \alpha$ for some universal 
constant $C > 0$.
Using again that $\mathrm{supp} \{ \nu_{\bar{y}}\} \subset \bar{\Omega} - \bar{y}$, applying~\eqref{xin}
and the arguments given in~\cite{BT, Topp} we conclude that $\Theta_2$ is uniformly away from the origin and vanishes as $\epsilon \to 0$. 
This last fact and (M2) lead us to the following estimate 
\begin{equation*}
\int_{\Theta_2} \nu_{\bar{y}}(dz) = o_\epsilon(1),
\end{equation*}
and therefore we arrive at
\begin{equation*}
\I_2^\delta \leq C \alpha + R o_\epsilon(1).
\end{equation*}

Using this,~\eqref{I_delta} and~\eqref{A} into~\eqref{testinglinearization} we can write
\begin{equation*}
\eta + \lambda_R M - o_\alpha(1) - o_\epsilon(1) - \epsilon^{-1}o_\delta(1) \leq o_\gamma(1).
\end{equation*}

Recalling that $\lambda_R$ is nonnegative, we arrive to a contradiction with $\eta > 0$ by 
taking $\delta \to 0, \epsilon \to 0, \gamma \to 0$ and $\alpha \to 0$.
\qed

As it is usual in the viscosity theory, comparison principle lead us to the following well-posedness result for~\eqref{eqparabolic}.
\begin{prop}\label{teoexistence}
Assume hypotheses of Proposition \ref{comparison} hold. 
Then, there exists a unique $u \in C(\bar{Q})$, viscosity solution to problem~\eqref{eqparabolic}.  
For each $T > 0$, this solution satisfies
\begin{equation*}\label{ubound}
|u(x, t)| \leq (||H(\cdot, 0, 0)||_{L^\infty(\bar{\Omega})} + ||f||_{L^\infty(\bar{Q}_T)}) \ t 
+ ||\varphi||_{L^\infty(\partial Q_T)} + ||u_0||_{L^\infty(\bar{\Omega})}, 
\end{equation*}
for all $(x,t) \in \bar{Q}_T$.
Moreover, if there exists $\lambda_0 > 0$ such that, for all $R > 0$, $\lambda_R \geq \lambda_0$,
and $f, \varphi$ are uniformly bounded,
then $u$ is uniformly bounded in $\bar{Q}$, with
\begin{equation*}
|u(x, t)| \leq \lambda_0^{-1} (||H(\cdot, 0, 0)||_{L^\infty(\bar{\Omega})} + ||f||_{L^\infty(\bar{Q})}) 
+ ||\varphi||_{L^\infty(\partial Q)} + ||u_0||_{L^\infty(\bar{\Omega})},
\end{equation*}
for all $(x,t) \in \bar{Q}$.
\end{prop}

The existence of a solution in the case $u_0, \varphi$ are smooth follows from the application of Perron's 
method for discontinuous solutions with generalized boundary conditions and comparison principle, see~\cite{usersguide, Ishii2, DaLio}.
The existence and uniqueness for general continuous initial and boundary data is obtained by approximation 
through smooth data and viscosity stability.
The $L^\infty$-estimates for the solution are easily obtained by taking sub and supersolutions for the problem with the 
form $(x,t) \mapsto \pm (C_1 t + C_2)$, with $C_1 \geq 0$ and $C_2 > 0$ suitable constants depending on the data.

\medskip

Following the same lines of the proof of Proposition~\ref{comparison} we can get the corresponding results for the state-constraint
evolution problem and for the stationary equation. Both results are presented next.
\begin{prop}\label{comparisons-c}
Under the assumptions of Proposition~\ref{comparison}, comparison principle holds for bounded viscosity sub and supersolutions
of the state-constraint problem
\begin{eqnarray}\label{eqcomparisons-c}
\left \{ \begin{array}{rll} \partial_t u - \I(u, x) + H(x,u,Du) & = f \quad & \mbox{in} \ Q_T, \\
\partial_t u - \I(u, x) + H(x,u,Du) & \geq f \quad & \mbox{in} \ \partial Q_T, \\
u & = u_0 \quad & \mbox{on} \ \bar{\Omega} \times \{ 0 \}. \end{array} \right .
\end{eqnarray}
\end{prop}

\begin{prop}\label{teostationary}
Let $\I$ be as in~\eqref{operator} defined through a family of measures $\{ \nu_x \}$ satisfying (M0)-(M3) with the same $\sigma <1$. Let $\varphi \in C(\partial \Omega)$, $H \in C(\bar{\Omega} \times \R \times \R^n)$ satisfying (H1) and (H2)
with $m>\sigma$. We consider the  problem
\begin{eqnarray}\label{inftyeq}
\left \{ \begin{array}{rcll} - \mathcal{I}(u) + H(x, u, Du) &=& 0 \quad & in \ \Omega, \\ 
u &=& \varphi \quad & in \ \partial \Omega, \end{array} \right . 
\end{eqnarray}

Further, assume one of the following additional hypotheses on $H$ holds

\medskip
\noindent
(i) There exists $\lambda_0 > 0$ such that, for all $R > 0$, we have $\lambda_R \geq \lambda_0$.
\medskip

\noindent
(ii) For each $x \in \bar{\Omega}$, the function $(r, p) \mapsto H(x,r, p)$ is convex 
and problem~\eqref{inftyeq} has a bounded strict subsolution.

\medskip

If $u, v$ are bounded, respective sub and supersolution to the problem~\eqref{inftyeq},
then $u \leq v$ in $\Omega$. Moreover, if we define
\begin{equation}\label{tildeustat}
\tilde{u}(x,t) = \left \{ \begin{array}{ll} u(x) \quad & \mbox{if} \ x \in \Omega,  \\
\limsup \limits_{y \to x, y \in \Omega} u(y) \quad & \mbox{if} \ x \in \partial \Omega,
\end{array} \right . 
\end{equation}
then $\tilde{u} \leq v$ in $\bar{\Omega}$.
In the above setting, there exists a unique solution $u \in C(\bar{\Omega})$ for~\eqref{inftyeq}.
\end{prop}

The following is a well-known consequence of comparison principle, see~\cite{Tchamba}.
\begin{lema}\label{Lipschitztime}
Assume that the hypotheses of Proposition~\ref{comparison} hold. Assume further that $\varphi \in C(\partial Q), f \in C(\bar{Q})$ 
are bounded and $u_0 \in C^1(\bar{\Omega})$. 
Then, the unique viscosity solution $u \in C(\bar{Q})$ 
to problem~\eqref{eqparabolic} satisfies
\begin{equation*}
|u(x,t) - u(x, s)| \leq C |t - s|, \quad \mbox{for all} \ s, t > 0; \ x \in \bar{\Omega},
\end{equation*}
where $C > 0$ depends only on $||u_0||_{C^1(\bar \Omega)}$ and the data.
\end{lema}

Combining the above lemma and Theorem~\ref{teolip} lead us to the following
\begin{cor}\label{corLiptime}
Assume that the hypotheses of Lemma~\ref{Lipschitztime} hold. Then, there exists a constant $C > 0$ just depending 
on the data and $||u_0||_{C^1}$ such that, the unique viscosity solution $u$ to~\eqref{eqparabolic} satisfies 
\begin{equation*}
|u(x,s) - u(y,t)| \leq C(|x - y| + |s - t|) \quad \mbox{for all} \ (x,s), (y,t) \in \bar Q.
\end{equation*}
\end{cor}


\section{The Ergodic Problem.}\label{sec:ep}

In this section we are interested in the existence and uniqueness of a constant $c \in \R$ for which the 
state-constraint problem
\begin{equation}\label{ergodic}
\left \{ \begin{array}{rll} -\I(u) + H_0(x, Du) &=  c \quad & \mbox{in} \ \Omega, \\ 
-\I(u) + H_0(x, Du) & \geq  c \quad & \mbox{on} \ \partial \Omega, \end{array} \right . 
\end{equation}
has a bounded viscosity solution. Here and below $H_0$ is a continuous function which satisfies (H2) with $m > 1$. Such a problem is known in the literature as the \textsl{ergodic problem}.
The key questions are the uniqueness of $c$ as well as the structure of the solutions $u$ to this problem.
For definition and further concepts associated to state-constraint problems we refer to~\cite{Fleming-Soner}.


\subsection{Solvability}
Existence of $c$ comes as a consequence of comparison and regularity results given in the previous sections.
\begin{prop}\label{ergodicexistence}
Let $\I$ as in~\eqref{operator} satisfying (M0)-(M3) with the same $\sigma <1$ and $H_0 \in C(\bar{\Omega} \times \R^N)$ 
satisfying (H2) with $m > 1$.
Then, there exists a unique constant $c \in \R$ for which the problem~\eqref{ergodic} has a viscosity solution 
in $C(\bar{\Omega})$, which is Lipschitz continuous on $\bar{\Omega}$. 
\end{prop}

Following Tchamba~\cite{Tchamba}, the above proposition can be obtained through the next lemmas.

\begin{lema}\label{lemaergodic}
Let $\I$ as in~\eqref{operator} satisfying (M0)-(M3) with the same $\sigma <1$, and $H_0 \in C(\bar{\Omega} \times \R^N)$ satisfying (H2).
For each $\alpha \in (0,1)$, there exists a unique viscosity solution $u_\alpha \in C(\bar{\Omega})$ to the problem
\begin{equation}\label{eqlemaergodic}
\left \{ \begin{array}{rll} \alpha u -\I(u) + H_0(x,Du) &= 0 & \quad \mbox{in} \ \Omega, \\ 
\alpha u -\I(u) + H_0(x, Du) & \geq 0 & \quad \mbox{on} \ \partial \Omega, \end{array} \right .
\end{equation}
and such a solution satisfies the inequality
\begin{equation}\label{ulambdabound}
||\alpha u_\alpha||_{L^\infty(\bar{\Omega})} \leq \tilde{C}, 
\end{equation}
for $\tilde{C} > 0$ independent of $\alpha$.
\end{lema}

\noindent
{\bf \textit{Proof:}} For $R > 0$, we are going to introduce the intermediate problem
\begin{equation}\label{eqR}
\left \{ \begin{array}{rll} \alpha u -\I(u) + H_0(x,Du) &= 0 & \quad \mbox{in} \ \Omega, \\ 
u & = R & \quad \mbox{on} \ \partial \Omega, \end{array} \right .
\end{equation}

Since $\alpha > 0$, by Proposition~\ref{teostationary} there exists a unique solution $u_{\alpha, R} \in C(\bar{\Omega})$ to~\eqref{eqR}.
Since the constant function equal to $-\alpha^{-1} ||H_0(\cdot, 0)||_{L^\infty(\bar{\Omega})}$ is a subsolution for this problem,
we have $u_{\alpha, R} \geq -\alpha^{-1} ||H_0(\cdot, 0)||_{L^\infty(\bar{\Omega})}$. The idea is to provide an upper bound which is independent 
of $R$. For this, we recall $\delta_0 > 0$ is such that the function $x \mapsto \mathrm{dist}(x, \partial \Omega)$ 
is smooth in the set $\Omega_{\delta_0}$.
We are going to denote $d \in C^2(\Omega) \cap C(\bar{\Omega})$ a nonnegative function defined as 
$d(x) = \mathrm{dist}(x, \partial \Omega)$ if $x \in \bar{\Omega}_{\delta_0}$
and which is strictly positive (depending only on $\delta_0$) in $\Omega \setminus \Omega_{\delta_0}$.
Then, for $\beta \in (0,1)$ and by Lemma~\ref{lemaId} we see that
\begin{equation*}
\I(d^\beta, x) \leq C d^{\beta - \sigma}(x), \quad \mbox{for all} \ x \in \Omega.
\end{equation*}

For $x \in \Omega$ and by the coercivity of $H_0$, evaluating the function $-d^\beta$ into~\eqref{eqR} we see that
\begin{equation*}
\begin{split}
& \alpha (-d^\beta(x)) - \I(-d^\beta, x) + H_0(x, -\beta d^{\beta - 1}(x) Dd(x)) \\
\geq & - \mathrm{diam}(\Omega)^\beta - C d^{\beta - \sigma}(x) + c \beta^m d^{m(\beta - 1)}(x) - C,
\end{split}
\end{equation*}
where $C, c > 0$ are universal constants. Recalling the coercivity degree of $H_0$ given by $m > \sigma$ in (H2), 
fixing $0 < \beta < (m - \sigma)/(m - 1)$ if $m > 1$ and any $\beta > 0$ if $m \leq 1$, 
and taking $\bar{\delta} < \delta_0 /4$ very small in terms of the data, 
we conclude $-d^\beta$ is a supersolution to~\eqref{eqR} in the set $\Omega_{\bar{\delta}}$. 
Moreover, this function is a supersolution up to $\partial \Omega$ in the generalized sense because there is no smooth function 
touching $-d^\beta$ from below at the boundary. Additionally, there exists a constant $C_1 \geq 0$ (depending only on the data and $\bar{\delta}$)
such that, for all $x \in \Omega \setminus \Omega_{\bar{\delta}/2}$ we have
\begin{equation*}
\lambda (-d\beta(x)) - \I(-d^\beta, x) + H_0(x, -\beta d^{\beta - 1}(x) Dd(x)) \geq -C_1.
\end{equation*}

Then, considering the function
\begin{equation*}
\psi(x) = 2\alpha^{-1}\Big{(} ||H_0(\cdot, 0)||_{L^\infty(\bar{\Omega})} + C_1\Big{)} - d^\beta(x),
\end{equation*}
we see by the above discussion that $\psi$ is a viscosity supersolution (in the generalized boundary sense) to~\eqref{eqR} and therefore,
by strong comparison principle, for all $R > 0$ we get the estimate
\begin{equation*}\label{estimateulambdaR}
-\alpha^{-1} ||H_0(\cdot, 0)||_{L^\infty(\bar{\Omega})} \leq u_{\alpha, R}
 \leq 2\alpha^{-1}\Big{(} ||H_0(\cdot, 0)||_{L^\infty(\bar{\Omega})} + C_1\Big{)}.
\end{equation*}

Hence, defining $\tilde{C} = 2 (||H_0(\cdot, 0)||_{L^\infty(\bar{\Omega})} + C_1)$ and taking $R > \tilde{C}/\alpha$,
the solution $u_{\alpha, R}$ to~\eqref{eqR} cannot satisfy the boundary condition for supersolutions in the 
classical sense and then it satisfies the state-constraint problem~\eqref{eqlemaergodic}. By uniqueness, we 
conclude the result labeling $u_\alpha = u_{\alpha, R}$ for all $R > \tilde{C}/\alpha$. The estimate~\eqref{ulambdabound} for $u_\alpha$
is inherited from the corresponding inequality for $u_{\alpha, R}$.
\qed

\begin{lema}\label{lemaliperg}
Let $\I$ as in~\eqref{operator} satisfying (M0)-(M3) and $H_0 \in C(\bar{\Omega} \times \R^N)$ satisfying (H2).
Then, for each $\alpha \in (0,1)$, the solution $u_\alpha \in C(\bar{\Omega})$ 
of the problem~\eqref{eqlemaergodic} is Lipschitz continuous in $\bar{\Omega}$ with Lipschitz constant depending 
on the data and $\mathrm{osc}_{\bar{\Omega}}(u_\alpha)$. 

Moreover, if $m > 1$ in (H2) the Lipschitz constant depends only on the data and not 
on $\alpha$ nor $||u_\alpha||_{L^\infty(\bar{\Omega})}$.
\end{lema}

The proof comes as a combination of~\eqref{ulambdabound} and Theorem~\ref{teolip}.
The additional property for superlinear Hamiltonians comes from Corollary~\ref{oscbound}.

\medskip
\noindent
{\bf \textit{Proof of Proposition~\ref{ergodicexistence}}:}
Fix $x^* \in \Omega$ and denoting $u_\alpha$
the unique solution to~\eqref{eqlemaergodic} given in Lemma~\ref{lemaergodic}, we define
\begin{equation*}
v_\alpha(x) = u_\alpha(x) - u_\alpha(x^*), \quad x \in \bar{\Omega},
\end{equation*}
which is a viscosity solution to
\begin{equation*}
\left \{ \begin{array}{rll} \alpha v_\alpha -\I(v_\alpha) + H_0(x,Dv_\alpha) &= -\alpha u_\alpha(x^*) & \quad \mbox{in} \ \Omega, \\ 
\alpha v_\alpha -\I(v_\alpha) + H_0(x,Dv_\alpha) &\geq -\alpha u_\alpha(x^*) & \quad \mbox{on} \ \partial \Omega, \end{array} \right .
\end{equation*}

By~\eqref{ulambdabound} we see that $\alpha u_\alpha(x^*)$ is uniformly bounded as $\alpha \to 0$. Thus, by Corollary~\ref{oscbound} and
Lemma~\ref{lemaliperg}, the family $\{ v_\alpha \}_{\alpha \in (0,1)}$ is uniformly bounded and equi-Lipschitz. 
Then, letting $\alpha \to 0$, classical stability results  
provide us the existence of a pair $(v, c) \in \mathrm{Lip}(\bar{\Omega}) \times \R$, viscosity solution to~\eqref{ergodic}.

Concerning the uniqueness of $c$, we consider $(v_1, c_1), (v_2, c_2) \in \mathrm{Lip}(\bar{\Omega}) \times \R$
solving the ergodic problem.
It is direct to see that for each $i=1,2$, the function $w_i(x,t) := v_i(x) - c_i t$ solves the parabolic 
state-constraint problem
\begin{eqnarray*}
\left \{ \begin{array}{rll} \partial_t w_i - \I(w, x) + H_0(x,Dw_i) & = 0 \quad & \mbox{in} \ Q_T, \\
\partial_t w_i - \I(w_i, x) + H_0(x,Dw_i) & \geq 0 \quad & \mbox{in} \ \partial Q_T, \\
w_i & = v_i \quad & \mbox{on} \ \bar{\Omega} \times \{ 0 \}. \end{array} \right  .
\end{eqnarray*}

Applying comparison principle given in Proposition~\ref{comparisons-c}, we see that
\begin{equation*}
w_1(x,t) \leq w_2(x,t) + ||v_1 - v_2||_{L^\infty(\bar{\Omega})} \quad \mbox{for all} \ (x,t) \in \bar{Q}.
\end{equation*}

Thus, $c_2 - c_1 \leq 2||v_1 - v_2||_{L^\infty(\bar{\Omega})}/t$ and we arrive at $c_2 \leq c_1$ making $t \to \infty$.
Exchanging the roles of $c_1$ and $c_2$ we conclude the uniqueness of $c$.
\qed

\subsection{Strong Maximum Principle.}
The aim of this subsection is to provide a version of the Strong Maximum Principle
for the following evolutive counterpart of~\eqref{ergodic}
\begin{equation}\label{parergodic}
\left \{ \begin{array}{rcll} \partial_t u - \mathcal{I}(u) + H_0(x, Du) &=& c \quad & \hbox{in} \ Q, \\ 
\partial_t u - \mathcal{I}(u) + H_0(x, Du) &\geq& c \quad & \hbox{on} \ \partial Q. \end{array} \right . 
\end{equation}

This Strong Maximum Principle is inspired by arguments from Coville~\cite{coville} (see also Ciomaga~\cite{ciomaga})
and plays a key role in the arguments to come. In fact, this is not only used to provide valuable information about 
the solutions of~\eqref{ergodic} (see Proposition~\ref{teoergodic} below), but it is a cornerstone in the proof of the 
large time behavior of our Cauchy-Dirichlet problem at several levels.

Recalling $H_0 \in C(\bar{\Omega} \times \R^N)$ trivially satisfies (H1), Proposition~\ref{comparisons-c}
provides us comparison principle for~\eqref{parergodic} and therefore we can get
\begin{lema}\label{lemakappa}
Assume (M0)-(M3) holds with the same $\sigma <1$ and $H_0 \in C(\bar{\Omega} \times \R^N)$ satisfies (H2).
Let $u$ u.s.c. in $\bar{Q}$, $v$ l.s.c. in $\bar{Q}$ with $u, v \in L^\infty(\bar{Q}_T)$ for all $T > 0$ 
be respective sub and supersolutions to problem~\eqref{parergodic}. 
For $t \in [0,+\infty)$, define 
\begin{equation}\label{kappa}
\kappa(t) = \sup \limits_{x \in \bar{\Omega}} \{ u(x, t) - v(x, t)\}. 
\end{equation}

Then, for all $0 \leq s \leq t$, we have $\kappa(t) \leq \kappa(s)$.
\end{lema}

\medskip

Next we introduce some notation: let $\{ \nu_x \}_x$ in the definition of $\I$ and 
for $x \in \R^N$ we define inductively
\begin{equation*}
X_0(x) = \{ x \}, \quad X_{n + 1}(x) = \bigcup_{\xi \in X_n(x)} \{ \xi + \mathrm{supp} \{ \nu_\xi \}) \}, \quad \mbox{for} \ n \in \N,
\end{equation*}
and
\begin{equation}
\mathcal{X}(x) = \overline{\bigcup_{n \in \N} X_n}. 
\end{equation}

The Strong Maximum Principle presented here relies in the nonlocality of the operator under the ``iterative covering property'' 
which is close to the ideas of Coville \cite{coville} but it has to be combined here with different arguments. 
This property is established through the condition
\begin{equation}\label{supportnu}
\mathcal{X}(x) = \bar{\Omega}, \quad \mbox{for all} \ x \in \bar{\Omega}.
\end{equation}

Notice this condition is satisfied by our two main examples~\eqref{censoredLaplacian} and~\eqref{regional}.


\begin{prop}\label{SMP}\textsc{(Strong Maximum Principle)}
Let $\I$ as in~\eqref{operator}, where $\{ \nu_x \}$ satisfies (M0)-(M3) and~\eqref{supportnu}, and 
$H_0 \in C(\bar{\Omega} \times \R^N)$ satisfying (H2).

\medskip
\noindent
\textsl{1.- Strong Maximum Principle - Parabolic Version:}
Let $u, v$ be respective bounded, viscosity sub and supersolution to~\eqref{parergodic}, with $u = \tilde{u}$ as in~\eqref{tildeu}.
Let $\kappa$ as in~\eqref{kappa} and assume there exists $t_0 > 0$ satisfying 
\begin{equation*}
\kappa(t_0) = \sup \limits_{t \geq 0} \{ \kappa(t) \}.
\end{equation*}

Then, the function $u - v$ is constant in $\bar{\Omega} \times [0, t_0]$. Moreover, we have
\begin{equation*}
(u - v)(x,t) = \kappa(0), \quad \mbox{for all} \ (x,t) \in \bar{\Omega} \times [0, t_0].
\end{equation*}

\medskip
\noindent
\textsl{2.- Strong Maximum Principle - Stationary Version:}
Let $u, v$ be respective bounded viscosity sub and supersolution to~\eqref{ergodic}, with $u = \tilde{u}$ as in~\eqref{tildeustat}. 
Then, $u - v$ is constant in $\bar{\Omega}$.
\end{prop}

\noindent
{\bf \textit{Proof:}} We focus in the parabolic version, and we start defining $T = t_0 + 1$. 
As we did it in the proof of comparison principle Proposition~\ref{comparison}, we may assume $u$ 
is Lipschitz continuous (in space and time) in $\bar{Q}_T$ by replacing $u$ by its sup-convolution $u^\gamma$. We avoid the direct
use of $u^\gamma$ for simplicity.

Notice that by Lemma~\ref{lemakappa}, we see that $\kappa(t) = \kappa(0)$ for all $t \in [0, t_0]$.
Fix $\tau \in (0, t_0)$, denote
$$
\mathcal{M}_\tau = \{ x \in \bar{\Omega} : (u - v)(x, \tau) = \kappa(\tau)\},
$$ 
and choose a point $x_\tau \in \mathcal{M}_\tau$. 
For simplicity, we also assume that $x_\tau \in \Omega$, otherwise
we apply Soner type penalization procedure used in the proof of the comparison principle instead of the function $\Phi$ below.
In that case, the control of the integral terms can be made in the same way as in the proof of the comparison principle.

Consider $\epsilon, \alpha > 0$ and define the function
\begin{equation*}
\Phi(x,y,s,t) = u(x,s) - v(y,t) - \phi(x,y,s,t),
\end{equation*}
with $\phi(x,y,s,t) := \epsilon^{-2} (|x - y|^2 + (s -t)^2) + \alpha((s - \tau)^2 + |x - x_\tau|^2)$. 
This function attains its maximum $(\bar{x}, \bar{y}, \bar{s}, \bar{t})$
in $\bar{Q}_T \times \bar{Q}_T$, and using the inequality 
$$
\Phi(\bar{x}, \bar{y}, \bar{s}, \bar{t}) \geq \Phi(x_0, x_0, \tau, \tau),
$$
we see that, keeping $\alpha > 0$,
\begin{equation}\label{SMPeps}
\bar{s}, \bar{t} \to \tau, \quad \bar{x}, \bar{y} \to x_\tau, \quad \mbox{and} \quad u(\bar{x}, \bar{s}) \to u(x_\tau, \tau), \quad 
v(\bar{y}, \bar{t}) \to v(x_\tau, \tau),
\end{equation}
as $\epsilon \to 0$. Moreover, by the Lipschitz continuity of $u$, $\bar{p} := 2 \epsilon^{-2} (\bar{x} - \bar{y})$ is uniformly 
bounded as $\epsilon \to 0$. Thus,  properly using $\phi$ as a test function for $u$ at $(\bar{x}, \bar{s})$ and for $v$ at $(\bar{y}, \bar{t})$
and substracting the corresponding viscosity inequalities, for each $\delta > 0$ we can write
\begin{equation}\label{SMPtesting}
2\alpha(\bar{s} - \tau) - (\epsilon^{-1} + \alpha)o_\delta(1) - \I^\delta + \mathcal{A} \leq 0,
\end{equation}
where
\begin{equation*}
\begin{split}
\I^\delta & := \I[B_\delta^c](u(\cdot, \bar{s}), \bar{x}) - \I[B_\delta^c](v(\cdot, \bar{t}), \bar{y}) \\
\mathcal{A} & := H_0(\bar{x}, \bar{p} + \alpha(\bar{x} - x_\tau)) - H_0(\bar{y}, \bar{p}). 
\end{split}
\end{equation*}

Note that by the continuity of $H_0$, using~\eqref{SMPeps} and the boundedness of $\bar{p}$, we readily conclude that
\begin{equation*}
\mathcal{A} = o_\epsilon(1) + o_\alpha(1). 
\end{equation*}

Now we deal with the nonlocal terms. Using (M3), the Lipschitz continuity of $u$ and~\eqref{SMPeps}, we can write
\begin{equation*}
\begin{split}
\I^\delta \leq & \ L o_\epsilon(1)
+ \int \limits_{B_\delta^c \cap (\Omega - \bar{x})} (u(\bar{x} + z, \bar{s}) - v(\bar{y} + z, \bar{t}) - (u(\bar{x}, \bar{s}) - v(\bar{y}, \bar{t})) \nu_{\bar{y}}(dz) \\
& \ - \int \limits_{B_\delta^c \setminus (\Omega - \bar{x})} (v(\bar{y} + z, \bar{t}) - v(\bar{y}, \bar{t}))\nu_{\bar{y}}(dz) \\
=: & \ L o_\epsilon(1) + \I_1^\delta - \I_2^\delta, 
\end{split}
\end{equation*}
where $L > 0$ is the Lipschitz constant for $u.$
But using that $x_\tau \in \Omega$,
by the second convergence in~\eqref{SMPeps} we see that $B_{\delta'} \subset (\Omega - \bar{x})$ for each $0 < \delta' \leq d(x_\tau)/2$.
By (M3) we have $\mathrm{supp} \{ \nu_{\bar{y}} \} \subset (\Omega - \bar{y})$ and using again the second convergence in~\eqref{SMPeps}
we obtain $|\mathrm{supp} \{ \nu_{\bar{y}} \} \setminus (\Omega - \bar{x})| \to 0$ as $\epsilon \to 0$. This
together with the boundedness of $v$ and (M1) allows us to write, for all $0 < \delta < \delta'$
\begin{equation*}
\I_2^\delta = o_\epsilon(1), 
\end{equation*}
meanwhile, using that $(\bar{x}, \bar{y}, \bar{s}, \bar{t})$ is maximum for $\Phi$ and (M2) lead us to
\begin{equation*}
\I_1^\delta \leq \alpha o_{\delta'}(1) 
+ \int \limits_{B_{\delta'}^c \cap (\Omega - \bar{x})} (u(\bar{x} + z, \bar{s}) - v(\bar{y} + z, \bar{t}) - (u(\bar{x}, \bar{s}) - v(\bar{y}, \bar{t})) \nu_{\bar{y}}(dz),
\end{equation*}

Replacing the above estimates for $\mathcal{A}$ and $\I^\delta$ in~\eqref{SMPtesting},
we make $\delta \to 0$ to get rid of the term $\epsilon^{-1}o_\delta(1)$. Recalling~\eqref{SMPeps},
by the boundedness assumptions for $u$ and $v$ and the continuity of the measure given by (M3), keeping $\delta' > 0$ 
and making $\epsilon \to 0$ and finally $\alpha \to 0$, Dominated Convergence Theorem leads us to
\begin{equation*}
- \int \limits_{B_{\delta'}^c \cap (\Omega - x_\tau)} (u(x_\tau + z, \tau) - v(x_\tau + z, \tau) 
- (u(x_\tau, \tau) - v(x_\tau, \tau)) \nu_{x_\tau}(dz) \leq 0,
\end{equation*}
but since $x_\tau \in \mathcal{M}_\tau$ and by the upper semicontinuity of $u - v$, we arrive at
\begin{equation*}
\int \limits_{B_{\delta'}^c \cap (\Omega - x_\tau)} (u(x_\tau + z, \tau) - v(x_\tau + z, \tau) 
- (u(x_\tau, \tau) - v(x_\tau, \tau)) \nu_{x_\tau}(dz) = 0. 
\end{equation*}

Thus, since $\delta' > 0$ is arbitrary, we conclude that $u - v \equiv \kappa(\tau)$ in $X_1(x_\tau)$ and proceeding 
inductively as above, we obtain 
\begin{equation*}
(u - v)(x, \tau) = \kappa(\tau) \quad \mbox{for each} \ x \in \bigcup_{n \in \N} X_n(x_\tau), 
\end{equation*}
from which we obtain $(u - v)(\cdot, \tau) = \kappa(\tau)$ in $\bar{\Omega}$ by applying~\eqref{supportnu}. The result 
for $\tau = 0$ and $\tau = t_0$ can be easily obtained by Lemma~\ref{lemakappa} and the upper semicontinuity of $u- v$.
\qed


As a consequence of the Strong Maximum Principle, we obtain the following
\begin{prop}\label{teoergodic}
Let $\I$ as in~\eqref{operator} satisfying (M0)-(M3) and the iterative covering property~\eqref{supportnu}, let 
$H_0 \in C(\bar{\Omega} \times \R^N)$ satisfying (H2) with $m > 1$ and $c$ the unique ergodic constant given 
in Proposition~\ref{ergodicexistence}.
Then, the solution to~\eqref{ergodic} is unique up to an additive constant.
\end{prop}

\noindent
{\bf \textit{Proof:}} Consider 
$(v_i, c) \in \mathrm{Lip}(\bar{\Omega}) \times \R$, $i =1,2$ be two
solutions to~\eqref{ergodic}. Then, since $v_1, v_2$ are time-independent, they can be seen as solution of \eqref{parergodic} and we may cast $v_1$ as a subsolution (in $\Omega$) to the problem for which $v_2$ is a supersolution. Moreover 
$$\kappa(t) = \sup \limits_{x \in \bar{\Omega}} \{ v_1(x) - v_2(x)\},$$
does not depend on $t$. 
Thus, we conclude $v_2 = v_1 + C$ for some $C \in \R$ by the Strong Maximum Principle (Proposition~\ref{SMP}).
\qed


\section{Large Time Behavior.}\label{sec:ltb}

For simplicity, in this section we concentrate on the parabolic problem
\begin{equation}\label{parabolic2}\tag{CP'}
\left \{ \begin{array}{rcll} \partial_t u + \lambda u - \mathcal{I}(u) + H_0(x, Du) &=& 0 \quad & \hbox{in} \ Q, \\ 
u &=& \varphi \quad & \hbox{on} \ \partial Q, \\
u &=& u_0 \quad & \hbox{on} \ \bar{\Omega} \times \{ 0 \}, \end{array} \right . 
\end{equation}
where $\lambda \geq 0$. In the rest of this section, we always assume $\I$ is as in~\eqref{operator} and 
satisfies (M0), (M1)-(M2) with the same $\sigma \in (0,1)$, and (M3). We also assume
$H_0 \in C(\bar{\Omega} \times \R^N)$ satisfies (H2) with $m > \sigma$, 
and $u_0 \in C(\bar{\Omega})$, $\varphi \in C(\partial \Omega)$ (hence, $\varphi$ is time independent) 
satisfy the compatibility condition~\eqref{compatibility}.
In this setting, problem~\eqref{parabolic2} can be uniquely solved in $C(\bar{Q})$.

For simplicity, assuming an initial data $u_0$ in~\eqref{parabolic2} the solution to~\eqref{parabolic2}
is Lipschitz continuous in $\bar{Q}$ with Lipschitz constant depending on the data and the $||u_0||_{C^1(\bar{\Omega})}$, 
see Corollary~\ref{corLiptime}. 
We get the large time behavior for the general case $u_0 \in C(\bar{\Omega})$ by approximation through smooth initial data, 
see~\cite{Tchamba} for a complete exposition of these arguments.

We start with the study of parabolic problems for which we have steady state asymptotic behavior.
\begin{teo}\label{teosteady}
Assume one of the following conditions hold

\medskip
\noindent
(i) $\lambda > 0$.

\medskip
\noindent
(ii) $\lambda = 0$, the ergodic constant $c$ associated to the ergodic problem~\eqref{ergodic} is negative, and the function
$p \mapsto H_0(x, p)$ is convex for all $x \in \bar{\Omega}$.

\medskip

In both cases, there exists a unique $C(\bar{\Omega})$-viscosity solution to the problem
\begin{equation}\label{S}\tag{$\mathrm{S}_\lambda$}
\left \{ \begin{array}{rcll} \lambda u - \mathcal{I}(u) + H_0(x, Du) &=& 0 \quad & in \ \Omega, \\ 
u &=& \varphi \quad & on \ \partial \Omega, \end{array} \right .
\end{equation}
and therefore, the unique solution to~\eqref{parabolic2} converges uniformly on $\bar{\Omega}$ 
as $t \to \infty$ to the unique viscosity solution to $(\mathrm{S}_\lambda)$.
\end{teo}

\noindent
{\bf \textit{Proof:}} The solvability of $(\mathrm{S}_\lambda)$ in the case $\lambda > 0$ comes as a consequence of the case $(i)$ 
in Proposition~\ref{teostationary}. In case $(ii)$, we notice that a solution 
to~\eqref{ergodic} is a strict subsolution to problem ($S_0$) and we fall in case $(ii)$ of Proposition~\ref{teostationary}, 
from which the solvability of ($S_0$) holds.

We note that in both cases $(i)$ and $(ii)$, the solution $u$ to~\eqref{parabolic2} is bounded in $\bar{Q}$. In fact, this can 
be easily seen in case $(i)$ through Proposition~\ref{teoexistence}. For case $(ii)$, we see that the function
\begin{equation*}
(x,t) \mapsto  u_\infty(x) - ||u_\infty - \varphi||_{L^\infty(\partial \Omega)} - ||u_\infty - u_0||_{L^\infty(\bar{\Omega})}, 
\end{equation*}
where $u_\infty$ is a solution to the ergodic problem~\eqref{ergodic},
is a visosity subsolution to the problem satisfied by $u$, and by comparison principle we conclude
\begin{equation*}
u_\infty - ||u_\infty - \varphi||_{L^\infty(\partial \Omega)} - ||u_\infty - u_0||_{L^\infty(\bar{\Omega})} 
\leq u \quad \mbox{on} \ \bar{Q}.
\end{equation*}

To find an upper bound for $u$, we fix 
$\bar{x} \in \Omega$ and denote $K = 2\mathrm{diam}(\Omega)$. By (H2) there exists a constant $C_K > 0$ depending 
on the data and $K$ such that the function
\begin{equation*}
 (x,t) \mapsto C_K(1 - K^{-1}|x - \bar{x}|) + ||\varphi||_{L^\infty(\partial \Omega)} + ||u_0||_{L^\infty}(\bar{\Omega}),
\end{equation*}
is a viscosity supersolution to~\eqref{parabolic2} and from this we see that
\begin{equation*}
u \leq C_K + ||\varphi||_{L^\infty(\partial \Omega)} + ||u_0||_{L^\infty}(\bar{\Omega}) \quad \mbox{on} \ \bar{Q}.
\end{equation*}

From this point, we argue simulteneously for cases $(i)$ and $(ii)$. The above analysis implies $u$ is uniformly bounded on $\bar{Q}$.
Then, the functions
\begin{equation*}
\bar{u}(x,t) = \limsup_{\epsilon \to 0, z \to x} u(z, t/\epsilon), \quad  
\underline{u}(x,t) = \liminf_{\epsilon \to 0, z \to x} u(z, t/\epsilon),
\end{equation*}
are well-defined in $\bar{Q}$. Naturally $\underline{u} \leq \bar{u}$ in $Q \cup \bar{\Omega} \times \{ 0 \}$.
Besides, for each $t > 0$, applying half-relaxed limits method~\cite{BP1, BP2} we see that the functions 
\begin{equation*}
x \mapsto \bar{u}(x, t) \quad \mbox{and} \quad x \mapsto \underline{u}(x,t)
\end{equation*}
are respective viscosity sub and supersolution to the problem~\eqref{S}. Thus, by comparison we obtain 
$\tilde{\bar{u}} \leq \underline{u}$ and then they coincide in $\Omega$. 

We claim that $\tilde{\bar{u}} = \bar{u}$ on $\partial \Omega$ and postpone its justification until the end of the 
convergence proof. If this happens, then the function 
\begin{equation*}
x \mapsto U(x, t) := \bar{u}(x, t) = \underline{u}(x, t) = \lim \limits_{z \to x, \epsilon \to 0} u(z, t/\epsilon), 
\quad x \in \bar{\Omega}
\end{equation*}
is in $C(\bar{\Omega})$ and it is a viscosity solution to~\eqref{S}, which is the unique one by Proposition~\ref{teostationary}.
From this, we easily conclude the uniform convergence on $\bar{\Omega}$ of the solution of the parabolic 
problem~\eqref{parabolic2} to the stationary problem~\eqref{S} as $t \to \infty$.

Now we deal with the claim. Let $x_0 \in \partial \Omega$. 
Note that by definition, for each $\eta > 0$ small, there exist $y_\eta, z_\eta \in \Omega$ 
and $\epsilon_\eta > 0$ satisfying $|x_0 - y_\eta|, |x_0 - z_\eta| \leq \eta$ and such that
\begin{equation*}
\bar{u}(x_0, t) - \tilde{\bar{u}}(x_0, t) \leq \eta + u(y_\eta, t/\epsilon_\eta) - u(z_\eta, t/\epsilon_\eta).
\end{equation*}

But we have that $u$ is uniformly Lipschitz in $\bar{Q}$ (see Lemma~\ref{Lipschitztime}) and then
\begin{equation*}
\bar{u}(x_0, t) - \tilde{\bar{u}}(x_0, t) \leq \eta + C||u_0||_{C^1(\bar{\Omega})} \eta,
\end{equation*}
for some $C > 0$ depending only on the data.
A similar lower bound can be stated and since $\eta$ is arbitrary, we conclude the equality.
\qed

\medskip

Now we address the ergodic large time behavior.
\begin{teo}\label{LTB}
Assume that $\lambda = 0$, $\I$ satisfies the iterative covering property \eqref{supportnu} and that 
$H_0$ satisfies (H2) with $m > 1$.
Let $u \in C(\bar{Q})$ be the unique solution to~\eqref{parabolic2}
and $(u_\infty, c) \in \mathrm{Lip}(\bar{\Omega}) \times \R$ be an ergodic solution to~\eqref{ergodic}. Then

\medskip
\noindent
\textsc{1.- Case I.} If $c > 0$, then $(S_0)$ has no bounded viscosity solution. 
Moreover, the function $(x,t) \mapsto u(x,t) + ct$ is uniformly bounded in $\bar{Q}$ and
\begin{equation*}
u + ct \to u_\infty + K \quad \mbox{uniformly in} \ \bar{\Omega}, \quad \mbox{as} \ t \to +\infty, 
\end{equation*}
for some constant $K \in \R$ depending on $H, u_0, \varphi$ and $c$.

\medskip
\noindent
\textsc{2.- Case II.} If $c = 0$, then any bounded solution to $(\mathrm{S}_0)$ has the form $u_\infty + K$ 
for some constant $K \in \R$ such that $u_\infty + K \leq \varphi$ on $\partial \Omega$. Moreover, 
\begin{equation*}
u \to u_\infty + K \quad \mbox{uniformly in} \ \bar{\Omega}, \quad \mbox{as} \ t \to +\infty, 
\end{equation*}
for some constant $K \in \R$ depending on $H, u_0, \varphi$ and $c$, and such that $u_\infty + K \leq \varphi$ on $\partial \Omega$.
\end{teo}

\noindent
{\bf \textit{Proof:}} \textsl{1.- Case I:} 
The nonexistence of a bounded viscosity solution for ($S_0$) can be proved by contradiction. If $\phi$ is
a bounded viscosity solution to ($S_0$), since $c > 0$ we have that $\phi$ is a strict subsolution to the equation~\eqref{ergodic} 
in $\Omega$. The function $u_\infty - 2(||u_\infty||_\infty + ||\phi||_\infty) - 1$ is a solution for this last problem, and therefore, by 
comparison principle we arrive at $\phi \leq u_\infty - 2(||u_\infty||_\infty + ||\phi||_\infty) - 1$ in $\bar{\Omega}$, which is a contradiction.

For the convergence, we start noting that the function $v(x,t) = u(x,t) + ct$ is a viscosity solution to the problem
\begin{equation*}
\left \{ \begin{array}{rcll} \partial_t v - \mathcal{I}(v) + H_0(x, Dv) &=& c \quad & \hbox{in} \ Q, \\ 
v &=& \varphi + ct \quad & \hbox{on} \ \partial Q, \\
v &=& u_0 \quad & \hbox{on} \ \bar{\Omega} \times \{ 0 \}. \end{array} \right . 
\end{equation*}

Now, considering the functions
\begin{equation*}
\psi_{\pm}(x,t) = u_\infty(x,t) \pm (||u_\infty||_{L^\infty(\bar{\Omega})} + ||\varphi||_{L^\infty(\partial \Omega)} + ||u_0||_{L^\infty(\bar{\Omega})}), 
\end{equation*}
we note that $\psi_-$ is a viscosity subsolution for the above problem satisfying the boundary condition in the classical sense because $c > 0$, 
meanwhile $\psi_+$ is a viscosity supersolution for the same problem satisfying the boundary condition in the generalized sense. Thus, we
obtain that $v$ is uniformly bounded in $\bar{Q}$ (by a constant depending on the data and $||u_\infty||_{L^\infty(\bar{\Omega})}$).

Now, consider $v$ as a subsolution to the equation
\begin{equation}\label{eqLBT}
\partial_t u - \mathcal{I}(u) + H_0(x, Du) = c \quad  \hbox{in} \ Q, 
\end{equation}
and $u_\infty$ as a supersolution for the state-constraint problem 
associated to this equation. Defining 
\begin{equation}\label{kappaLBT}
\kappa(t) = \max \limits_{x \in \bar{\Omega}} \{ v(x,t) - u_\infty(x) \}, \quad t \geq 0, 
\end{equation}
we can apply Lemma~\ref{lemakappa} to conclude that $\kappa$ is decreasing in $t$. 
Since it is also bounded, we conclude that $\kappa(t) \to \bar{\kappa} \in \R$, as $t \to \infty$.

Considering 
$\{ v_k \}_{k}$ the sequence of functions in $\mathrm{Lip}(\bar{Q})$ given by $v_k(x,t) = v(x, t + t_k)$, we see that 
\begin{equation*}
|v_k(x,t) - v_k(x', t')| \leq C|x - x'| + (C + c) |t - t'|, \quad \mbox{for each} \ (x,t), (x',t') \in \bar{Q},
\end{equation*}
where $C > 0$ depends on the data and the smoothness of $u_0$. Then, $\{ v_k \}_k$ is equicontinuous and bounded, and by Arzela-Ascoli Theorem,
up to subsequences, $v_k \to \bar{v} \in \mathrm{Lip}(\bar{Q})$, 
uniformly in $\bar{Q}_T$ for each $T > 0$. Stability results for viscosity solutions imply that the limit
function $\bar{v}$ satisfies~\eqref{eqLBT}. Evaluating~\eqref{kappaLBT} in $t + t_k$ and taking limit
as $k \to \infty$ we arrive at
\begin{equation*}
\bar{\kappa} = \max \limits_{x \in \bar{\Omega}} \{ \bar{v}(x,t) - u_\infty(x) \} \quad \mbox{for each } \ t \geq 0,
\end{equation*}
and since $u_\infty$ is a solution for the state-constraint problem associated to~\eqref{eqLBT}, by the application of the Strong Maximum Principle
given by Proposition~\ref{SMP}, we conclude that
\begin{equation*}
\bar{v}(x,t) = u_\infty(x) + \bar{\kappa}, \quad \mbox{for all} \ (x,t) \in \bar{Q}. 
\end{equation*}

This, together with the definition of $\bar{v}$ and the fact that $\kappa$ is nonincreasing allows us to conclude the result.

\medskip
\medskip
\noindent
\textsl{2.- Case II:} Let $\phi \in C(\bar{\Omega})$ be a bounded solution to ($S_0$). Since $u_\infty$ is a 
solution to~\eqref{ergodic}, 
it can be regarded as a viscosity supersolution to ($S_0$) with generalized boundary condition. Thus, the stationary version of the 
Strong Maximum Principle leads us to the existence of a constant $K \in \R$ such that $\phi = u_\infty + K$. 
But the stationary version of Lemma~\ref{lemasubboundary} implies that $\phi \leq \varphi$ on $\partial \Omega$, from which we conclude the result.
The convergence proof follows the same lines of the previous case.
\qed

\bigskip

\noindent
{\bf Acknowledgements:} G.B. is partially supported by the ANR (Agence Nationale de
la Recherche) through ANR WKBHJ (ANR-12-BS01-0020). 
E.T. was partially supported by CONICYT, Grants Capital Humano Avanzado and Cotutela en el 
Extranjero.


\end{document}